\documentclass{article}
\usepackage{arxiv}

\usepackage[utf8]{inputenc} 
\usepackage[T1]{fontenc}    
\usepackage{hyperref}       
\usepackage{url}            
\usepackage{booktabs}       
\usepackage{amsfonts}       
\usepackage{nicefrac}       
\usepackage{microtype}      
\usepackage{lipsum}

\usepackage{graphicx}

\usepackage{lineno,hyperref}
\usepackage{amssymb}
\usepackage{amsmath,color}
\usepackage{amsfonts}
\usepackage{dsfont}

\usepackage{graphicx}
\usepackage{amsmath}%

\newtheorem{theorem}{Theorem}

\newtheorem{corollary}[theorem]{Corollary}

\newtheorem{example}[theorem]{Example}

\newtheorem{lemma}[theorem]{Lemma}

\newtheorem{proposition}[theorem]{Proposition}
\newtheorem{remark}[theorem]{Remark}

\newcommand{\beq} {\begin{eqnarray*}}
\newcommand{\eeq} {\end{eqnarray*}}

\newcommand{\noi} {\noindent}

\DeclareMathOperator*{\argmin}{Argmin}

\def \R{\mathbb{R}}

\def \F{\mathbb{F}}
\def \G{\mathbb{G}}

\title{Weak convergence of empirical Wasserstein type distances}
\author{P. Berthet \thanks{
              Institut de Math\'{e}matiques de Toulouse ; UMR 5219
 Universit\'{e} de Toulouse ; CNRS
 UPS IMT, F-31062 Toulouse Cedex 9,
              \texttt{philippe.berthet@math.univ-toulouse.fr} }         
           \And
          J.-C. Fort \thanks{
              MAP5 ; UMR 8145 
 Universit\'{e} Paris Descartes ; CNRS ;
 45 rue des Saints p\`eres, F-75006,
\texttt{jean-claude.fort@parisdescartes.fr} }}

\begin{document}





\maketitle

\begin{abstract}
We estimate contrasts $\int _0 ^1 \rho(F^{-1}(u)-G^{-1}(u))du$ between two continuous distributions $F$ and $G$ on $\mathbb R$ such that the set $\{F=G\}$ is a finite union of intervals, possibly empty or $\mathbb{R}$. The non-negative convex cost function $\rho$ is not necessarily symmetric and the sample may come from any joint distribution $H$ on $\mathbb{R}^2$ with marginals $F$ and $G$ having light enough tails with respect to $\rho$. The rates of weak convergence and the limiting distributions are derived in a wide class of situations including the classical Wasserstein distances $W_1$ and $W_2$. The new phenomenon we describe in the case $F=G$ involves the behavior of $\rho$ near $0$, which we assume to be regularly varying with index ranging from $1$ to $2$ and to satisfy a key relation with the behavior of $\rho$ near $\infty$ through the common tails. Rates are then also regularly varying with powers ranging from $1/2$ to $1$ also affecting the limiting distribution, in addition to $H$.

\texttt{Central limit theorems, Generalized
Wasserstein distances, Empirical processes, Strong approximation, Dependent samples, Non-parametric statistics, Goodness-of-fit tests}.\\
{62G30, 62G20, 60F05, 60F17}
\end{abstract}

\section{Introduction}\label{sec:intro}
\subsection{Motivation}
In \cite{BFK17} we addressed the problem of estimating the distance between two asymptotically well separated and continuous distributions on the real line $\mathbb{R}$, with respect to a large class of generalized Wasserstein costs. The framework was the same as in \cite{Freitag} and is very simple. A sequence of independent and indentically distributed ($i.i.d.$) random variables ($r.v.$) taking values in $\mathbb R^2$ is available. The marginals have distinct continuous cumulative distribution function ($c.d.f.$) $F$ and $G$. For instance, each couple may result from simultaneous experiments. We estimated contrasts $\int_0^1 c(F^{-1}(u),G^{-1}(u))du$ between $F$ and $G$ by the natural and easily computed non-parametric plug-in estimator $\int_0^1 c(\F_n^{-1}(u),\G_n^{-1}(u))du$.  Here $F^{-1}$ is the generalized inverse of $F$, $\F_n$ is the empirical $c.d.f.$, and $c$ is a non-negative cost. The almost sure ($a.s.$) consistency of this estimator being easily established under minimal assumptions we mainly developed a sharp method of proof of the Central Limit Theorem (CLT) assuming that the tails of $F$ and $G$ are distinct enough and compatible with the cost $c$. The most original contribution in \cite{BFK17} was to investigate rather deeply the latter relationship in the untrimmed case and for dependent samples. This showed that the problem can not be reduced to the study of each marginal $\int_0^1 c(\F_n^{-1}(u),F^{-1}(u))du$ and instead requires crossed assumptions on tails, costs and densities beyond moments. However the special case of the distance $W_1$ was not captured, asymptotically non-symmetric costs or asymptotically too close marginals were not allowed, the case $F=G$ and the one marginal case were not considered.

In the present paper -- the first version of this preprint is \cite{BF18} -- the general setting remains exactly the same, but we investigate the most important situations for statistical applications, among which the goodness-of-fit hypothesis $F=G$, the alternative hypothesis where $F\ne G$ on $\mathbb R$ and may have arbitrarily close tails, and the intermediate hypothesis where the two situations $F=G$ and $F\ne G$ are encountered, but alternate along a finite number of intervals. The distance $W_1$ and non-symmetric costs are now allowed provided that they are regularly varying at both sides of $0$. We focus on the new difficulties, however we often refer to \cite{BFK17} to borrow some long arguments and apply already developed tools. New assumptions arise that again illustrate how delicate tail integrals of transforms of empirical quantile functions can be for heavy-tailed distributions.

The method of proofs relies on a careful subdivisions of the integrals and events, and a joint approximation of the quantile processes $\sqrt n(\F_n^{-1}(u)-F^{-1}(u))$, $u\in(0,1)$ by properly scaled Brownian bridges on an appropriate sub-interval. As a matter of fact, it is not possible to directly apply a functional delta-method since the Hadamard differentiability of $F \to F^{-1}$ can not be extended to encompass distribution with densities arbitrarily close to $0$ and in particular with unbounded supports. Moreover the Brownian approximation - weak or strong - of the quantile processes suffer many problems near $0$ and $1$ due to extreme values. Lastly, the general costs we use - even the simple Wasserstein costs - make the problem more difficult to handle and shows up to be determinant for both rates and limits in the case $F=G$.

Let us mention related results in the framework of univariate probability distributions. The commonly used $p$-Wasserstein distance $W_p(F,G)$ is
\begin{equation}\label{Frechet}
W_p^p(F,G)=\int_0^1|F^{-1}(u)-G^{-1}(u)|^p du.
\end{equation}
Many authors were interested in the convergence of $W_p^p(\mathbb{F}_n,F)$, see $e.g.$ the survey paper \cite{Ledoux17} or \cite{Delbarrio99,Delbarrio05,Delbarrio11}. Up to our knowledge there are only two recent works studying the convergence of $W_2^2(\mathbb{F}_n,\mathbb{G}_n)$ \cite{dBL,Munk}, for independent samples. The results of \cite{dBL} are valid in any finite dimension with the drawback that the estimator is not explicit from the data and the centering in the central limit theorem (CLT) is the biased $\mathbb EW_2^2(\mathbb F_n,\mathbb G_n)$ rather than $W_2^2(F,G)$ itself, moreover the limiting variance has no closed form expression and seems not easy to estimate. In \cite{Munk} the estimator is the same as our's, howewer only discrete distributions and $W_2$ distance are considered. Notice also that in the early work \cite{CzadoMunk} a trimmed version of the Mallows distance $W_2^2(\mathbb F_n,\mathbb G_n)$ is studied, however under an implicit assumption on the level of trimming which has to hold in probability. Moreover in the case of dependent samples, a trimmed version of $W_2^2(\mathbb{F}_n,\mathbb{G}_n)$ is studied in \cite{Freitag}.

We investigate below a larger class of convex costs, even larger than in \cite{BFK17}. The samples are possibly not independent, and the conditions relating the tails of $F$ and $G$ to the cost function $c$ are easily checked. Combined to our technique of proof they allow to control the critical parts of the untrimmed integrals in a weaker sense than in probability, hence our explicit sufficient conditions are lighter than the above mentionned implicit ones. We obtain a general CLT for $W_c(\F_n,\G_n)$ when $F=G$ are continuous, thus providing a new class of goodness-of-fit and comparison tests with exact rates and non-degenerate limits. In order to evaluate the power of these tests we study the weak convergence under many alternatives $F\ne G$ among which the case where $F=G$ on large intervals.

\subsection{Setting}
The $p$-Wasserstein distance between two $c.d.f.$ $F$ and $G$ on $\mathbb R$ is defined by
\begin{equation}\label{Wass}
W^p_p(F,G)= \min_{X\sim F,Y\sim G}\mathbb E|X-Y|^p
\end{equation}
where $X\sim F,Y\sim G$ means that X and Y are joint real $r.v.$ having $c.d.f.$ $F$ and $G$ respectively. The minimum in (\ref{Wass}) is (\ref{Frechet}). To any non negative function $c(x,y)$ from $\mathbb R^2$ to $\mathbb R$ let associate the Wasserstein type cost
\begin{equation}\label{Wassc}
W_c(F,G)= \min_{X\sim F,Y\sim G}\mathbb E c(X,Y).
\end{equation}
We are interested in triplets $(c,F,G)$ such that $W_c(F,G)$ is finite and can be estimated by using an explicit CLT. To guaranty that an analogue of \eqref{Frechet} exists we consider cost functions defining a negative measure on $\mathbb{R}^2$, hence satisfying
\begin{equation}\label{P}
c(x',y')-c(x',y)-c(x,y')+c(x,y)\leqslant 0, \quad x\leqslant x', y\leqslant y'.
\end{equation}
If $c$ satisfies (\ref{P}) then for any functions $a$ and $b$, $a(x)+b(y)+c(x,y)$ satisfies (\ref{P}). In particular $c(x,y)=-xy$ and $(x-y)^2=x^2+y^2-2xy$ satisfy (\ref{P}). More generally if $\rho$ is a convex real function then $c(x,y)=\rho(x-y)$ satisfies (\ref{P}). Two important special cases are the symmetric power functions $|x-y|^p$, $p\geqslant 1$, associated to $W_p$ and the non-symmetric contrast functions $c(x,y)=(x-y) (\alpha-{\bf 1}_{x-y<0})$ associated to the $\alpha ^{th}$ quantile, $0< \alpha <1$. The following result yields the minimum in (\ref{Wassc}) in a closed form analogous to \eqref{Frechet}.

\begin{theorem}[Cambanis, Simon, Stout \cite{Cambanis76}]\label{theo:cam} If $c$ satisfies (\ref{P}) 
then 
$$ W_c(F,G)=\int_0^1 c(F^{-1}(u),G^{-1}(u))du.
$$
\end{theorem}

Let $\displaystyle (X_i,Y_i)_{1\leqslant i\leqslant n}$ be an $i.i.d.$ sample of a random vector with joint $c.d.f.$ $H$ on $\mathbb R^2$ and marginal $c.d.f.$ $F$ and $G$ on $\mathbb R$. Write $\F_n$ and $\G_n$ the random empirical $c.d.f.$ built from the two marginal samples. Thus $\F_n$ and $\G_n$ are not independent in general. Consider a cost function $c$ satisfying (\ref{P}). Let $X_{(i)}$ (resp. $Y_{(i)}$) denote the $i^{th}$ order statistic of the sample $\displaystyle (X_i)_{1\leqslant i\leqslant n}$ (resp. $\displaystyle (Y_i)_{1\leqslant i\leqslant n}$), i.e. $X_{(1)}\leqslant \ldots\leqslant X_{(n)}$. By Theorem \ref{theo:cam}, the non-parametric statistic
\begin{equation}\label{WC}
W_c(\F_n,\G_n)=\frac{1}{n}\sum_{i=1}^n c(X_{(i)},Y_{(i)})
\end{equation}
is a natural estimator of $W_c(F,G)$. Now, the $c(X_{(i)},Y_{(i)})$'s being neither independent nor with identical distributions the statistic (\ref{WC}) is not classical - such as $i.i.d.$ mean, L-statistic, U-statistic etc. Notice also that $W_c(F,G)$ does not depend on the generally unknown $H$ whereas the $r.v.$ $W_c(\F_n,\G_n)$ strongly depends on $H$ through its distribution. In \cite{BFK17} we established the CLT
\[
\sqrt{n}\left(  W_{c}(\mathbb{F}_{n},\mathbb{G}_{n})-W_{c}(F,G)\right)
\rightarrow_{weak}\mathcal{N}\left(  0,\sigma^{2}\right)
\]
whenever the tails of $F$ and $G$ differ from at least $\tau >0$ and $c(x,y)$ is asymptotically $\rho(x-y)$ with $\rho$ non-negative, symmetric, convex. The influence of $H$ only appeared in the limiting variance $\sigma^{2}=\sigma^{2}(H,c)$ together with $c$. The sufficient conditions relating explicitly $c$, $F$ and $G$ were designed to  carefully control the extremes, define sharply the truncation level and approximate the underlying joint quantile processes. We now intend to complete the picture by extending this CLT to other important cases, in particular $\tau =0$ and non symmetric costs $\rho$. 

\subsection{Overview}

Hereafter we consider a cost $c(x,y)=\rho_c(x-y)$ where $\rho_c$ is a non-negative real convex function such that $\rho_c(0)=0$, and is not assumed to be symmetric. In the spirit of \cite{BFK17} we separate out three sets of assumptions, labeled $(FG)$, $(C)$ and $(CFG)$ respectively.

First, $(FG)$ concerns the regularity and tails of $F$ and $G$, and especially their density-quantile function. Conditions $(FG)$ are satisfied by distributions having regular tails, among which all classical probability distributions.

Second, $(C)$ restricts the rate of increase at infinity of $\rho_c$ and the regular variation at $0$ of $\rho_c$, without even assuming differentiability at $0$. Conditions $(C)$ encompass a large class of Wasserstein type costs $c$ and the distance $W_1$ is now allowed, together with non-symmetric variants of Wasserstein distances $W_p^p$, $p\geqslant 1$, possibly with slowly varying factors -- a non trivial extension -- or exponential and over-exponential costs. 

The conditions $(FG)$ and $(C)$ are thus designed to separately select a class of probability distributions and admissible costs.

The third set $(CFG)$ aims at mixing  the requirements on $c$, $F$ and $G$ making them compatible. We distinguish between $(CFG_E)$, $(CFG_D)$ and $(CFG_{ED})$ depending on the situations $\{F=G\}=\mathbb R$ or $\{F\ne G\}=\mathbb R$ or $\{F=G\}\ne\mathbb R$ and $\{F\ne G\}\ne\mathbb R$, respectively. The joint distribution $H$ of the couples is not restricted and again only affects the limiting distributions. In order to exhibit an exact rate of convergence it shows up that the tail constraints on $F$ and $G$ that naturaly depend on $\rho$ at $\infty$ also strongly depend on the exact regular variation of $\rho_c$ at $0$ whenever $F=G$ in tails, that is the key requirement of $(CFG_E)$ and $(CFG_{ED})$.

When dealing with empirical Wasserstein type integrals, to adapt the functional delta method one would need to truncate and then to assume a convergence in probability of the extremal parts. This would be a restriction excluding many distributions $F$ and $G$, depending on where the integral is non-adaptively trimmed. Moreover, proving the validity of the assumed convergence of the truncated parts would require variants of Steps 1, 2, 3 of our proofs. In contrast, $(CFG_{E})$ and $(CFG_{D})$ explicitly relate the tails to the cost in such a way that the implicit truncation levels can be defined appropriately.

Before entering the mathematical details of these assumptions let us present two consequences of our results. The regular variation of tails is in the sense of (i) in Section \ref{regul} below and $\rightarrow_{weak}$ denotes the convergence in distribution. 
\begin{proposition}\label{PropWp}
Consider the Wasserstein distance $W_p^p$ for $1<p<2$. Assume that $F=G$ is two times differentiable, $\log F(x)$ and $\log (1-F(x))$ are regularly varying as $|x|\to \infty$, and $\displaystyle F(x)(1-F(x))\leqslant C |x|^{-(\frac{2(p+2)}{2-p}+\varepsilon)}$ for some $\varepsilon>0$, $C>0$ and all $|x|$ large enough. Then it holds$$
\displaystyle
n^{\frac{p}{2}} W_{p}^p(\mathbb{F}_{n},\mathbb{G}_{n})\rightarrow_{weak}\int_{0}%
^{1}\left\vert \mathbb{B}(u)\right\vert
^{p}du
,$$
where $\mathbb B$ is an explicit centered Gaussian process and the limiting $r.v.$ is positive and finite. 
\end{proposition}

The restriction $p<2$ is not surprising since when $X$ and $Y$ are Gaussian and the two samples are independent, the limiting random integral is $a.s.$ infinite. More precisely, in the case $p=2$ we establish the weak convergence of $\displaystyle n W_{2}^2(\mathbb{F}_{n},\mathbb{G}_{n})$ by requiring $F$ to be sub-Gaussian, as in \cite{Delbarrio05} for $\displaystyle n W_{2}^2(\mathbb{F}_{n},F)$.

In the case $p=1$ we get, with the same Gaussian process $\mathbb B$ as above, the following result, which seems new for $W_{1}(\mathbb{F}_{n},\mathbb{G}_{n})=||\mathbb{F}_{n}-\mathbb{G}_{n}||_1$.
\begin{proposition}\label{PropW1}
Assume that the set $\{F=G\}$ is a finite union of non empty intervals of $\R$, that $F,G$ are two times differentiable and that  $\log F(x)$, $\log G(x)$, $\log (1-F(x))$ and $\log (1-G(x))$ are regularly varying as $|x|\to \infty$. Let $r=2$ if $\{F=G\}$ is compact, and $r=6$ otherwise. Assume that $\displaystyle \max (F(x)(1-F(x)),G(x)(1-G(x)))\leqslant C|x|^{-(r+\varepsilon)}$ for some $\varepsilon>0$, $C>0$ and all $|x|$ large enough. Then it holds
\[
\sqrt{n}\left(  W_{1}(\mathbb{F}_{n},\mathbb{G}_{n})-W_{1}(F,G)\right)
\rightarrow_{weak}\int_{F^{-1}\ne G^{-1}}\mathbb{B}(u)du+\int_{F^{-1}= G^{-1}}\left\vert \mathbb{B}%
(u)\right\vert du
\]
and the limiting $r.v.$ is finite.
\end{proposition}

As can be seen in the two previous results this paper focuses on the probability distributions with infinite support. Nevertheless our results also hold for compactly supported probability distributions with derivable densities. At the end of Section \ref{sec:main} we provide simplified sufficient assumptions in the compactly supported case.

The paper is organized as follows. Assumptions are discussed in Section \ref{sec:hyp}. In Section  \ref{sec:main} we state our main results in the form of CLT's for 
$\displaystyle W_{c}(\F_{n},\G_{n})-W_{c}(F,G)$ at various rates. We propose a few perspectives for applications in Section \ref{conclusion}. All the results are proved in Section~\ref{sec:proof}.  

\section{Assumptions}\label{sec:hyp}

\subsection{Assumptions $(FG)$}
Consider a sequence $(X_{n},Y_{n})\in \R^{2}$ 
of independent random vectors having the same $c.d.f.$ $H$ as $(X,Y)$. The distribution $H$ may have a density or not. However we assume that the marginal $c.d.f.$'s $F$ of $X$ and $G$ of $Y$ have support $\mathbb{R}$ and positive densities
$f=F^{\prime}$ and $g=G^{\prime}$. Let $(E,D)$ be the partition of $(0,1)$ defined by
\begin{equation}\label{defED}
E=\left\{  u:F^{-1}(u)=G^{-1}(u)\right\}  ,\quad D=\left\{  u:F^{-1}(u)\neq
G^{-1}(u)\right\}.
\end{equation}
If $u$ shifts infinitely many times between $E$ and $D$ it becomes difficult
to control the stochastic integral $W_{c}(\mathbb{F}_{n},\mathbb{G}_{n})$. The
case where $\left\vert F^{-1}(u)-G^{-1}(u)\right\vert >\tau>0$ as
$u\rightarrow1$ and $u\rightarrow0$ has been treated in details in \cite{BFK17}.
 We allow the diagonal $\left\vert
F^{-1}(u)-G^{-1}(u)\right\vert \leqslant\tau$ and thus encompass the case
$E=\left(  0,1\right)  $ together with some tractable situations where
$E\neq\emptyset$ and $D\neq\emptyset$. Let assume that there exists a finite
integer $\kappa\geqslant2$ and $0=u_{0}<u_{1}<...<u_{\kappa}=1$ such that,
writing $A_{k}=(u_{k-1},u_{k})$,%
\[
(FG0)\quad F^{-1}(u_{k})=G^{-1}(u_{k})\text{ and }A_{k}\subset
E\text{ or }A_{k}\subset D\text{, for }k=1,...,\kappa.
\]
This covers three generic cases, namely $E=\left(  0,1\right)  $, $D=\left(
0,1\right) $ and when $D\neq\emptyset$, $E \neq\emptyset$ are finite unions of intervals. The exponential rate of decrease of the right and left tails of $F$ and
$G$ are defined to be, for $x\in\mathbb{R}_{+}$,%
\begin{align*}
\psi_{X}^{+}(x) &  =-\log\mathbb{P}(X>x),\quad\psi_{Y}^{+}(x)=-\log
\mathbb{P}(Y>x),\\
\psi_{X}^{-}(x) &  =-\log\mathbb{P}(X<-x),\quad\psi_{Y}^{-}(x)=-\log
\mathbb{P}(Y<-x).
\end{align*}
Only $\psi_{X}^{+}$ and $\psi_{Y}^{+}$ will be considered in subsequent proofs
where arguments given for the right hand tail $u\rightarrow 1$ in the
integrals $W_{c}(F,G)$ and $W_{c}(\mathbb{F}_{n},\mathbb{G}_{n})$ work similarly for the left hand tail $u\rightarrow 0$. Define the density quantile functions $h_{X}=f\circ F^{-1}$ and $h_{Y}=g\circ G^{-1}$ then assume 
\[%
\begin{array}
[c]{l}%
(FG1)\quad F,G\in\mathcal{C}_{2}(\mathbb{R)}\text{,}\quad f,g>0\text{ on
}\mathbb{R}.\medskip\\
(FG2)\quad\sup\limits_{0<u<1}\min(u,1-u)\left\vert \left(  \log
h(u)\right)  ^{\prime}\right\vert <+\infty\quad\text{for }h=h_{X},h_{Y}.\\
(FG3)\quad\sup\limits_{0<u<1}\dfrac{\min(u,1-u)}{\left(  \left\vert
\Gamma^{-1}(u)\right\vert +1\right)  h(u)}<+\infty\quad\text{for }%
(h,\Gamma)=(h_{X},F)\text{ or }(h_{Y},G).
\end{array}
\]
Observe that $(FG1)$ and $(FG2)$ are classical in the context of approximation of quantile processes -- see e.g. \cite{CH93}.

\begin{remark}
Rewriting $(FG2)$ and $(FG3)$ we get%
\begin{align*}
\sup_{x\in\mathbb{R}}\frac{\min(F(x),1-F(x))}{f(x)}\left(  \frac
{1}{\left\vert x\right\vert +1}+\frac{|f^{\prime}(x)|}{f(x)}\right)   &
<+\infty,\\
\sup_{x\in\mathbb{R}}\frac{\min(G(x),1-G(x))}{g(x)}\left(  \frac
{1}{\left\vert x\right\vert +1}+\frac{|g^{\prime}(x)|}{g(x)}\right)   &
<+\infty.
\end{align*}
In Proposition 5 of \cite{BFK17} we provided a simple sufficient condition for $(FG1)$, $(FG2)$, $(FG3)$ based on the regular variation of
$\psi_{X}^{\pm}$ and $\psi_{Y}^{\pm}$.  All classical tail distributions satisfy the conditions $(FG)$.
\end{remark}

\subsection{Notation for regularity}\label{regul}

To specify the allowed cost functions $c(x,y)$ the following
definitions are required. As usual for $k\in\mathbb{N}_{\ast}$ and $I\subset\mathbb{R}$
let $\mathcal{C}_{k}(I)$ denote the set of functions that are $k$ times
continuously differentiable on $I$ and $\mathcal{C}_{0}(I)$ the set of
continuous functions on $I$. In forthcoming assumptions and proofs we consider
functions defined either on $\left(  0,x_{0}\right)  $ or on $\left(
y_{0},+\infty\right)  $ for some $0<x_{0}<y_{0}$. We distinguish the two
domains by using a variable $x\rightarrow0$ and a variable $y\rightarrow
+\infty$. In \cite{BFK17}   only large values $y\in\left(
y_{0},+\infty\right)  $ played a role in terms of regular variation, so that
we keep the same setting in (i) below. Unexpectedly, it turns out that the two
domains interfere when $|F-G|$ is arbitrarily small, and we need (ii).$\smallskip$

\noindent\textbf{(i)} Regularity on $\left(  y_{0},+\infty\right)  $. Let
$\mathcal{M}_{2}(\left(  y_{0},+\infty\right)  )$ be the subset of functions
$l\in\mathcal{C}_{2}(\left(  y_{0},+\infty\right)  )$ such that $l^{\prime
\prime}$ is monotone on $\left(  y_{0},+\infty\right)  $. Write $RV(+\infty
,\gamma)$ the set of regularly varying functions at $+\infty$ with index
$\gamma\geqslant0$. If $\gamma=0$ we restrict ourselves to slowly varying
functions $L$ at $+\infty$ such that%
\begin{equation}
L^{\prime}(y)=\frac{\varepsilon(y)L(y)}{y},\quad
\lim_{y\rightarrow+\infty}\varepsilon(y)=0.\label{L'}%
\end{equation}
This weak restriction is explained at Section 6 of \cite{BFK17}.
In order to find distributions $F$ and $G$ compatible with the cost $c$ we further impose
\begin{equation}
L^{\prime}(y)\geqslant\frac{l_1}{y},\quad l_1\geqslant1,\quad y\geqslant
y_{0}.\label{L1}%
\end{equation}
For $\gamma=0$, introduce%
\[
RV_{2}(+\infty,0)=\left\{  L:L\in\mathcal{M}_{2}\left(  \left(  y_{0}%
,+\infty\right)  \right)  \text{\ such that (\ref{L'}), (\ref{L1})
hold}\right\}
\]
and for $\gamma>0,$%
\[
RV_{2}(+\infty,\gamma)=\left\{  l:l\in\mathcal{M}_{2}\left(  \left(
y_{0},+\infty\right)  \right)  ,l(y)=y^{\gamma}L(y)\text{\ such that
}L^{\prime}\text{ obeys (\ref{L'})}\right\}  .
\]
\noindent\textbf{(ii)} Regularity on $\left(  0,x_{0}\right)  $. We consider
positive slowly varying functions $L$ at $0$,%
\begin{equation}
\lim_{x\searrow0}\frac{L(\theta x)}{L(x)}=1\text{ for any }\theta>0.\label{L0}%
\end{equation}
For $b>1$ let introduce%
\[
RV_{2}(0,b)=\left\{  \rho:L\in\mathcal{C}_{2}\left(  \left(  0,x_{0}\right)
\right)  ,\rho(x)=x^{b}L(x)\text{ such that }L\text{ satisfies (\ref{L0}%
)}\right\}  .
\]
For $b=1$ let define%
\[
RV_{2}(0,1)=\left\{  \rho:L\in\mathcal{C}_{2}\left(  \left(  0,x_{0}\right)
\right)  ,\rho(x)=xL(x)\text{ such that }L\text{ satisfies (\ref{L0}),
(\ref{Lpi})}\right\}
\]
where we impose the following finite limit%
\begin{equation}
\lim_{x\searrow0}L(x)=L(0)\in\mathbb{R}_{+}.\label{Lpi}%
\end{equation}

\subsection{Assumptions $(C)$}

We consider costs such that, for some $0<x_{0}%
<y_{0}<+\infty$,%
\[%
\begin{array}
[c]{l}%
({C0})\quad c(z,z^{\prime})=\rho_{c}(z-z^{\prime})\geqslant0,\quad
z,z^{\prime}\in\mathbb{R},\quad c(0,0)=0,\quad\rho_{c}\text{ is convex.}%
\smallskip\\
({C1})\quad\rho_{c}(x)=\rho_{-}(-x)1_{x\leqslant0}+\rho_{+}(x)1_{x\geqslant
0},\quad x\in\mathbb{R}, \quad \rho_{\pm}\in \mathcal{C}_{2}((0,+\infty)).\smallskip\\
({C2})\quad\rho_{+}(x)=x^{b_{+}}L_{+}(x)>0,\quad 0 < x \leqslant x_{0},\quad\rho_{+}\in{RV}_{2}(0,b_{+}\mathbb{)},\quad b_{+}\geqslant1,\\
\quad\quad\quad\rho_{-}(x)=x^{b_{-}}L_{-}(x)>0,\quad 0 < x \leqslant x_{0},\quad\rho_{-}\in{RV}_{2}(0,b_{-}\mathbb{)},\quad b_{-}\geqslant
1.\smallskip\\
({C3})\quad\rho_{+}(y)=\exp(l_{+}(y)),\quad y \geqslant y_{0},\quad l_{+}%
\in{RV}_{2}(+\infty,\gamma_{+}\mathbb{)},\mathbb{\quad}\gamma_{+}\geqslant0,\\
\quad\quad\quad\rho_{-}(y)=\exp(l_{-}(y)),\quad y \geqslant y_{0},\quad
l_{-}\in{RV}_{2}(+\infty,\gamma_{-}\mathbb{)},\quad\gamma_{-}\geqslant0.
\end{array}
\]
Notice that $\rho_{\pm}(0)=0$ and $\rho_{\pm}$ are  positive, continuous, convex and
increasing on $\mathbb{R}_{+}$. Define $\rho\left(  x\right)  =\max(\rho_{+}(x),\rho_{-}(x))$ and $b=\min(b_{+},b_{-})$. For $0\leqslant
x\leqslant x_{0}$ it holds%
\begin{equation}
\rho (x)=x^{b}L(x),\quad L(x)=\left\{
\begin{array}
[c]{cc}%
L_{+}(x) & \text{if }b_{+}<b_{-},\\
L_{-}(x) & \text{if }b_{-}<b_{+},\\
\max(L_{+}(x),L_{-}(x)) & \text{if }b_{+}=b_{-}.
\end{array}
\right.
\label{L}%
\end{equation}
Further assume that%
\[
({C4})\quad\lim_{x\rightarrow0}\frac{\rho_{+}(x)}{\rho(x)}=\pi_{+},\quad
\lim_{x\rightarrow0}\frac{\rho_{-}(x)}{\rho(x)}\rightarrow\pi_{-},\quad\pi
_{+},\pi_{-}\in\left[  0,1\right]  .
\]
Typical costs satisfying the conditions $(C)$ are the following.
\begin{example}
Let $a=(a_{-},a_{+})$ be such that $a_{\pm}>0$ and $b=(b_{-},b_{+})$ be such that $b_{\pm
}\geqslant1$. Then
$$
c_{a,b}(z,z^{\prime})=a_{-}\left(  z^{\prime}-z\right)  ^{b_{-}}%
1_{z<z^{\prime}}+a_{+}\left(  z-z^{\prime}\right)  ^{b_{+}}1_{z^{\prime}%
<z}%
$$
satisfies $(C)$ with $\gamma_- = \gamma_+=0$ and $\varepsilon(y)=O(1/\log y)$. This includes the Wasserstein distance $W^p_p$, $p\geqslant 1 $, by taking $a=(1,1)$ and $b=(p,p)$. It is possible to define costs mixing the Wasserstein distance $W^p_p$, $p\geqslant 1 $ near $0$ and $W^q_q$, $q\geqslant 1 $
away from $0$. Note that de facto when $E$ is not compact we will restrict to $p<2$ near $0$ in order to include at least the Gaussian distributions in $(CFG_E)$ and $(CFG_{ED})$ below. For instance the cost $\rho(x)=|x|(1+|x|)$ is well suited for distributions with heavier tails than Gaussian. 
\end{example}

\subsection{Assumptions $(CFG)$}

The joint influence of $l_{\pm}$, $L_{\pm}$ and $b_{\pm}$ on the
allowed tails $F^{-1}$ and $G^{-1}$ is expressed as follows. Remind the sets $E$ and $D$ from (\ref{defED}). We need three
different assumptions, each corresponding to the generic cases $E=\left(
0,1\right)  $, $D=\left(  0,1\right) $ and when at least one interval is
included in $E$ and one in $D$.

Studying the case $E=(0,1)$ we worked out the following conditions $(CFG_E)$. They only deal with the behavior of $F,G,\rho_c$ at infinity but also involve the orders $b_{\pm}\geqslant 1$ of the local regular variation $(C2)$ near zero that indeed rule the CLT rate. The case $b_-=2$ or $b_+=2$, which is restricted to sub-Gaussian distributions, is treated separately at Theorem \ref{W2}.\smallskip

\noindent\textbf{Assumption }$(CFG_{E})$. Assume that $b_{-}<2$ and $b_{+}<2$. Assume that for some $\theta_{2}>0$ and
\begin{equation}
(l,\psi)\in\left\{  (l_{+},\psi_{X}^{+}),(l_{-},\psi_{X}^{+}),(l_{-},\psi
_{X}^{-}),(l_{+},\psi_{X}^{-})\right\}  \label{lpsy}%
\end{equation}
we have,\\ 
$(i)$ if $1<b<2$, for all $y>y_{0}$,
\begin{equation}\label{lrondpsy}
l\circ\psi^{-1}(y)\leqslant\left(  1-\frac{b}{2}\right)  y+\log L\left(
\exp(-y/2)\right)  -2\log\psi^{-1}(y)-\theta_{2}\log y,
\end{equation}
and,\\ 
$(ii)$ if $b=1$, for all $y>y_{0}$,
\begin{equation}\label{lrondpsyb1}
l\circ\psi^{-1}(y)\leqslant\frac{y}{2}-2\log\psi^{-1}(y)-\theta_{2}\log
y.
\end{equation}
$\smallskip$

From the study of the case $D=\left(  0,1\right)$ the conditions $(CFG_D)$ that comes out only deal with the behavior of $F,G,\rho_c$ at infinity and the CLT rate is standard. The special case where $\left\vert F^{-1}(u)-G^{-1}(u)\right\vert >\tau>0$ as
$u\rightarrow1$ and $u\rightarrow0$ under $(C2)$ with $b>1$ is already covered by \cite{BFK17}. In order to cover more cases we further impose (\ref{CFGD_integree}) and allow $b=1$. Therefore $(CFG_D)$ extends the condition $(CFG)$ in \cite{BFK17}.
$\smallskip$

\noindent\textbf{Assumption }$(CFG_{D})$. Let $\theta_{-},\theta_{+}$ be the parameter $\theta>1$ of condition $(CFG)$
in \cite{BFK17} for the left and right tails respectively. \\
$(i)$ For any
$(l,\psi)$ from (\ref{lpsy}) and $\theta=\theta_{+}$ if $l=l_{+}$ or
$\theta=\theta_{-}$ if $l=l_{-}$ we have\\
\begin{equation}
(\psi\circ l^{-1})^{\prime}(y)\geqslant2+\frac{2\theta}{y},\quad y>y_{0}.\label{CFGD_derivee}
\end{equation}
$(ii)$ 
If \[
\underset{u\rightarrow1}{\lim\inf}\left\vert F^{-1}(u)-G^{-1}(u)\right\vert
=0\quad\text{or}\quad\underset{u\rightarrow0}{\lim\inf}\left\vert
F^{-1}(u)-G^{-1}(u)\right\vert =0
\]
and for $(l,\psi)$=$(l_{+},\psi_{X}^{+}),(l_{-},\psi_{Y}^{+})$ or $(l,\psi)=(l_{-},\psi_{X}^{-}),(l_{+},\psi_{Y}^{-})$  respectively, assume that for some $\theta_{2}>0$ it holds
\begin{equation}
l\circ\psi^{-1}(y)\leqslant\frac{y}{2}-2\log\psi^{-1}(y)-\theta_{2}\log
y,\quad y>y_{0}.\label{CFGD_integree}
\end{equation}$\smallskip$

When $D\neq\emptyset$ and $E\neq\emptyset$, two situations arise. Firstly, if $E$ is compact in $(0,1)$, that is $(A_{1}\cup A_{\kappa})\subset D$ we only need ($CFG_D$).  Secondly if at least one among $A_{1}$ or $A_{\kappa}$ is included in $E$, which means that $F=G$ on an infinite interval, then we need to also impose ($CFG_E$) on the involved intervals.\\

\noindent\textbf{Assumption }$(CFG_{ED})$. Assume $(CFG_{D})$. If $A_{1}\subset E$ then assume $(CFG_{E})$ for $(l,\psi)$=$(l_{-},\psi_{X}^{-}),(l_{+},\psi_{X}^{-})$. If $A_{\kappa}\subset E$ then assume $(CFG_{E})$ for $(l,\psi)=(l_{-},\psi_{X}^{+}),(l_{+},\psi_{X}^{+})$.\smallskip

\begin{remark}
If $\gamma_{\pm}>0$ we have $\theta_{\pm}>2$ and, if $\gamma_{\pm}=0$ we have, as in \cite{BFK17},
\[
\theta^{\pm}>2-\underset{y\rightarrow+\infty}{\lim\inf}\,\frac{\log
(1/\varepsilon_{\pm}(y))}{\log l_{\pm}(y)}%
\]
where $\varepsilon_{\pm}(y)$ corresponds to the function $\varepsilon(y)$ of
(\ref{L'}) applied to $L(y)=l_{\pm}(y)$.
\end{remark}

\begin{remark}
As will be seen in the proofs, $(CFG_{E})$ and $(FG3)$ imply that we can find $b^{\prime}$ such that
$1\leqslant b<b^{\prime}<2$ and
\begin{equation}
\int_{0}^{1}\left(  \frac{\sqrt{u(1-u)}}{h_X(u)}\right)  ^b
du\leqslant\int_{0}^{1}\left(  \frac{\left\vert F^{-1}(u)\right\vert }%
{\sqrt{u(1-u)}}\right)  ^{b^{\prime}}du<+\infty\label{bobkovledoux}%
\end{equation}
which is a little stronger than the necessary condition that the left hand integral is finite. By using
$F^{-1}(u)=\psi^{-1}(\log(1/(1-u)))$, (\ref{lrondpsy}) also reads%
\[
(F^{-1}(u))^{2}\rho\left(  F^{-1}(u)\right)  \leqslant\frac{L\left(
\sqrt{1-u}\right)  }{(1-u)^{1-b/2}(\log(1/(1-u)))^{\theta_{2}}},\quad u>u_{0}.
\]
In particular, if $L(x)=1$ we deduce that $(FG)$ and $(CFG_{E})$ imply
\[
\mathbb{P}\left(  X>y\right)  \leqslant\left(  \frac{1}{y^{2}\rho(y)}\right)
^{2/(2-b)},\quad y>y_{0}.
\]
This induces the moment conditions of Propositions \ref{PropWp} and \ref{PropW1}.
\end{remark}

\begin{example}
For light tails of Weibull type, $\psi(y)=y^{w}$, $w>0$, (\ref{bobkovledoux})
is true and $(CFG_{E})$ requires that $l(y^{1/w})<Cy$ as $y\rightarrow+\infty$
and hence a cost of type $l(y)=y^{\gamma}$, $y>y_{0}$ and $l(x)=x^{b}$,
$x<x_{0},$ is allowed provided that $\gamma<w$ and $1\leqslant b<2$. For heavy
tailed distributions such as Pareto, $\psi(y)=p\log y$ with index $p>2$, the conditions
$(CFG_{E})$, $(CFG_{D})$ and $(CFG_{ED})$ induce more constraints. For
instance $(CFG_{E})$ applied with $\rho_c(x)=x^{b}$, $x<x_{0},$ and $l(y)=\alpha
\log y$, $y>y_{0},$ implies that $p>4/(2-b)$ and $1\leqslant\alpha
<p(1-b/2-2/p)$, hence the minimal requirement on $p$ is $p\geqslant6/(2-b)$.
Choosing $\rho_c(x)=x^{b}$ on $\mathbb{R}_{+}$ we have $\alpha=b$ and the last
constraint becomes $p>2(b+2)/(2-b)$.
\end{example}

\section{Statement of the results}\label{sec:main}

Consider the joint Gaussian process $\displaystyle\mathcal{G}=\left\{  \left(  \mathbb{B}^{X}(u),\mathbb{B}^{Y}(u)\right)
:u\in(0,1)\right\}$ with
\begin{align} \label{gauss} \quad\mathbb{B}^{X}(u)=\frac{B^{X}(u)}{h_{X}(u)}%
,\quad\mathbb{B}^{Y}(u)=\frac{B^{Y}(u)}{h_{Y}(u)},
\end{align}
where $(B^{X},B^{Y})$ are two standard Brownian bridges with covariance
\[
cov(B^{X}(u),B^{X}(v))=cov(B^{Y}(u),B^{Y}(v))=\min(u,v)-uv,\quad u,v\in\left(
0,1\right)  ,
\]
and cross covariance
\[
cov(B^{X}(u),B^{Y}(v))=H(F^{-1}(u),G^{-1}(v))-uv,\quad u,v\in\left(
0,1\right)  .
\]
The existence of $\mathcal{G}$ is proved in \cite{BFK17}. Let
$\mathbb{B}(u)=\mathbb{B}^{X}(u)-\mathbb{B}^{Y}(u)$, $u\in(0,1)$, that is the Gaussian process driving the limit distribution in Propositions \ref{PropWp} and \ref{PropW1} as well as in forthcoming results.\medskip

We are now ready to state our main results. Remind (\ref{L}) and set
\begin{equation}
v_{n}=\frac{1}{\rho\left(  1/\sqrt{n}\right)  }=\frac{n^{b/2}}{L\left(
1/\sqrt{n}\right)  }\label{vn}%
\end{equation}
hence, in our first statement we have $K\sqrt{n}\leqslant v_{n}=o(n)$ for some $K>0$. The constants $\pi_{-}$ and $\pi_{+}$ come from $(C4)$. Our first result concerns $F=G$.

\begin{theorem} \label{E=R}
Assume $(FG)$, $(C)$, $E=\left(  0,1\right)  $ and $(CFG_{E})$, in which case $1\leqslant b^-,b^+<2$. Then%
\[
v_{n}W_{c}(\mathbb{F}_{n},\mathbb{G}_{n})\rightarrow_{weak}\pi_{-}\int_{0}%
^{1}1_{\left\{  \mathbb{B}(u)<0\right\}  }\left\vert \mathbb{B}(u)\right\vert
^{b_{-}}du+\pi_{+}\int_{0}^{1}1_{\left\{  \mathbb{B}(u)>0\right\}  }\left\vert
\mathbb{B}(u)\right\vert ^{b_{+}}du
\]
and the limiting $r.v.$ is finite and, if $\mathbb{P}(X=Y)<1$, positive.
\end{theorem}

\begin{remark}
As shown in \cite{CHS93}, and since $\mathbb{B}^X$ is a centered Gaussian process,
$$ 
\mathbb{P} \left( \int_{0}^{1} \left\vert \mathbb{B}^X(u)\right\vert ^{b} du < +\infty \right) = 1 \textnormal{ is equivalent to } \int_{0}^{1}\left(  \frac{\sqrt{u(1-u)}}{h_X(u)}\right)  ^{b}du < +\infty .
$$
The latter bound being guaranteed by $(CFG_{E})$ and $(FG3)$, which imply (\ref{bobkovledoux}), the finiteness of the limiting $r.v.$ in Theorem \ref{E=R} follows.
\end{remark}

For light tails one can
handle the limiting case $b=2$ -- here stated with $b=b^+=b^-=2$ and $L(x)=1$ for $|x|<x_{0}$ for sake of simplicity.

\begin{theorem} \label{W2}
Assume that $E=(0,1)$, $(FG1)$, $(FG2)$ and%
\begin{equation}\label{hypW2}
\lim_{u\rightarrow0}\frac{u}{h_X(u)}=\lim_{u\rightarrow1}\frac{1-u}%
{h_X(u)}=0,\quad\int_{0}^{1}\frac{u(1-u)}{h_X^{2}(u)}du<+\infty.
\end{equation}
Moreover assume $(C0)$ with $\rho_{c}(x)=x^{2}$ for $\left\vert x\right\vert
\leqslant x_{0}$, hence $b=b^+=b^-=2$. Then
$$
nW_{c}(\mathbb{F}_{n},\mathbb{G}_{n})\rightarrow_{weak}\int_{0}^{1}%
\mathbb{B}(u)^{2}du.
$$
\end{theorem}

Notice that Theorem \ref{W2} includes the case $W_c=W_{2}^2$ and shows that the cost function only matters at $0$.
\begin{example}
For light tails of Weibull type it holds, for some $w>0$,
$$h_X(u)=w(1-u)\left(\log(1/(1-u))\right)^{1-1/w}$$
and $(1-u)/h_X^2(u)=1/w\left((1-u)\log(1/(1-u))\right)^{2(1-1/w)}$. The first condition in (\ref{hypW2}) is then satisfied for $w>1$ and the second for $w>2$, so that $w>2$ is required. This excludes Gaussian tails, as in Theorem 4.6 in \cite{Delbarrio05}.
\end{example}

\begin{remark}
Theorem \ref{W2} requires no assumption on the cost $\rho(y)$ as $y\rightarrow+\infty$. In particular, $(C3)$ may hold with any $\gamma_+,\gamma_-$. Since only sub-Gaussian tails are allowed by (\ref{hypW2}) the tail part of $W_{c}(\mathbb{F}_{n},\mathbb{G}_{n})$ indeed behaves the same as for compactly supported distributions. Namely, empirical extremes of both samples remain simultaneously stuck together very closely to their common deterministic counterpart $F^{-1}$ that increases very slowly.
\color{black}
\end{remark}

Our second main statement is an extension of the main theorem of \cite{BFK17} which now allows $F$ and $G$ to have arbitrarily close tails.

\begin{theorem}\label{D=R}
Assume $(FG)$, $(C)$, $D=\left(  0,1\right)  $ and $(CFG_{D})$. Then%
\[
\sqrt{n}\left(  W_{c}(\mathbb{F}_{n},\mathbb{G}_{n})-W_{c}(F,G)\right)
\rightarrow_{weak}\mathcal{N}\left(  0,\sigma^{2}\right)
\]
where%
\[
\sigma^{2}=\mathbb{E}\left(  \left(  \int_{0}^{1}\vert\rho_{c}^{\prime}%
(F^{-1}(u)-G^{-1}(u))\vert\mathbb{B}(u)du\right)  ^{2}\right)  <+\infty.
\]

\end{theorem}

\begin{remark}
The finiteness and a closed form expression for $\sigma^{2}=\sigma^{2}(c,H)$ have been
proved in \cite{BFK17}. We also refer to the latter paper for explicit
examples in the independent samples case.
\end{remark}

Our third result shows that if there exists a point, or equivalently an open interval, where $F\ne G$ then the rate is $\sqrt n$, whether $E\neq\emptyset$ or not.

\begin{theorem}\label{ED}
Assume $(FG)$, $(C)$, $D\neq\emptyset$ and $(CFG_{ED})$. If $1<b<2$ then%
\[
\sqrt{n}\left(  W_{c}(\mathbb{F}_{n},\mathbb{G}_{n})-W_{c}(F,G)\right)
\rightarrow_{weak}\mathcal{N}\left(  0,\sigma_{D}^{2}\right)
\]
where%
\[
\sigma_{D}^{2}=\mathbb{E}\left(  \left(  \int_{D}\vert\rho_{c}^{\prime}%
(F^{-1}(u)-G^{-1}(u))\vert\mathbb{B}(u)du\right)  ^{2}\right)  <+\infty.
\]
If $b=1$ then, for $L_{\pm}(0)$ from (\ref{Lpi}),
\begin{align*}
\sqrt{n}\left(  W_{c}(\mathbb{F}_{n},\mathbb{G}_{n})-W_{c}(F,G)\right)   &
\rightarrow_{weak}\int_{D}\vert\rho_{c}^{\prime}%
(F^{-1}(u)-G^{-1}(u))\vert\mathbb{B}%
(u)du\\
&  +1_{\left\{  b_{-}=1\right\}  }L_{-}(0)\int_{E}1_{\left\{  \mathbb{B}%
(u)<0\right\}  }\left\vert \mathbb{B}(u)\right\vert du\\
&  +1_{\left\{  b_{+}=1\right\}  }L_{+}(0)\int_{E}1_{\left\{  \mathbb{B}%
(u)>0\right\}  }\left\vert \mathbb{B}(u)\right\vert du.
\end{align*}

\end{theorem}

\begin{remark}
In the second part of Theorem \ref{ED} the first term in the limiting $r.v.$ has distribution $\mathcal{N}\left(  0,\sigma_{D}^{2}\right)$ and is correlated in an explicit way to the other two terms. Theorem \ref{ED} also shows that whenever $1<b<2$ Theorem \ref{D=R} remains true if F and G are not stochastically ordered but cross each other at a finite number of points, since this implies $\sigma^{2}_D=\sigma^{2}$.
\end{remark}

The next corollary concerns the $L_1$-distance $W_{1}(\mathbb{F}_{n},\mathbb{G}_{n})=\Vert \mathbb{F}_{n}-\mathbb{G}_{n} \Vert _{L_1} $. Remind that $c_{a,1}(z,z^{\prime})=a_{-}\left(  z^{\prime}-z\right)
1_{z<z^{\prime}}+a_{+}\left(  z-z^{\prime}\right)  1_{z^{\prime}<z}$.

\begin{corollary} \label{pinball}
Assume $(FG)$, $(C)$ and $(CFG_{ED})$. Then
\begin{align*}
& \sqrt{n}\left(  W_{c_{a,1}}(\mathbb{F}_{n},\mathbb{G}_{n})-W_{c_{a,1}%
}(F,G)\right)\\    & \rightarrow_{weak}\int_{D}\left(  a_{-}1_{\left\{
F^{-1}(u)<G^{-1}(u)\right\}  }+a_{+}1_{\left\{  F^{-1}(u)>G^{-1}(u)\right\}  }\right)
\mathbb{B}(u)du\\
& +\int_{E}\left(  a_{-}1_{\left\{  \mathbb{B}(u)<0\right\}  }+a_{+}%
1_{\left\{  \mathbb{B}(u)>0\right\}  }\right)  \left\vert \mathbb{B}%
(u)\right\vert du
\end{align*}
and, in particular for $a_{-}=a_{+}=1$,%
\[
\sqrt{n}\left(  W_{1}(\mathbb{F}_{n},\mathbb{G}_{n})-W_{1}(F,G)\right)
\rightarrow_{weak}\int_{D}\mathbb{B}(u)du+\int_{E}\left\vert \mathbb{B}%
(u)\right\vert du.
\]
\end{corollary}

It is easily seen that straightforward adaptations of the proof of Theorems \ref{E=R} to \ref{W2} leads to analog results for $\sqrt{n} \left(W_{c}(\mathbb{F}_{n},G)-W_{c}(F,G)\right)$ and $v_n W_{c}(\mathbb{F}_{n},F)$ by just replacing $\mathbb B(u)=\mathbb{B}^{X}(u)-\mathbb{B}^{Y}(u)$ with $\mathbb{B}^{X}(u)$. In particular we get the following corollary of Theorem \ref{E=R}.

\begin{corollary}\label{WpE=R}
Let $1\leqslant p< 2$. Assume that $F$ satisfies $(FG)$ and has tails lighter than a Pareto tail with index strictly larger than $2(p+2)/(2-p)$. Then%
\[
n^{p/2}W_{p}^p(\mathbb{F}_{n},F)\rightarrow_{weak}\int_{0}%
^{1}\left\vert \mathbb{B}^X(u)\right\vert
^{p}du
\]
and the limiting $r.v.$ is positive and finite.
\end{corollary}

We conclude this section by stating the counterpart of Theorem \ref{E=R} for compactly supported probability distributions. Other extensions to this case of the above results are likewise easy.

\begin{corollary}\label{compact} 
Assume $wlog$ that $F=G$ has support $[0,1]$ and is twice differentiable with positive derivative $f$ on $(0,1)$. Assume moreover $(FG2)$, $(FG3)$ and $(C)$ except $(C3)$ with $b^-<b'$ and $b^+<b'$ where $b'>1$ and
\begin{equation}\label{eqcompact}
\int_{0}^{1}\left(  \frac{\sqrt{u(1-u)}}{h_X(u)}\right)  ^{b^{\prime}%
}du<+\infty.\end{equation}
Then \[
v_{n}W_{c}(\mathbb{F}_{n},\mathbb{G}_{n})\rightarrow_{weak}\pi_{-}\int_{0}%
^{1}1_{\left\{  \mathbb{B}(u)<0\right\}  }\left\vert \mathbb{B}(u)\right\vert
^{b_{-}}du+\pi_{+}\int_{0}^{1}1_{\left\{  \mathbb{B}(u)>0\right\}  }\left\vert
\mathbb{B}(u)\right\vert ^{b_{+}}du
\]
and the limiting $r.v.$ is finite and, if $\mathbb{P}(X=Y)<1$, positive.
\end{corollary}

This extends Theorem 19 of \cite{BFK17} to the case $F=G$ and reduces $(CFG_E)$ to the integrability assumption with no restriction on  $b$, since the influence of the cost is limited to its behaviour near $0$.

\begin{example} The Beta distribution with parameters $\alpha > 0$ and $\beta > 0$, has density $f(x)=\mathcal{B} (\alpha,\beta) \ x^{\alpha-1}(1-x)^{\beta-1}$ on $(0,1)$. Clearly $(FG2)$ and $(FG3)$ are satisfied, and since (\ref{eqcompact}) is true for any $b'> 1$, the previous result applies for any $b^-,b^+\geqslant1$.\\
This is not always the case. For instance, consider a $c.d.f.$ $F$ on $(0,1)$ equal to $\displaystyle e^{-1/|\log x|^w}$, $w>0$ near $0$ -- and symmetrically near $1$. Then it satisfies $(FG2)$ and $(FG3)$ but only satisfies (\ref{eqcompact}) for $b'\leqslant 2$. Hence the previous result  applies for $1\leqslant b^-$, $b^+<2$.


\end{example}

\section{Applications}\label{conclusion}

\subsection{Comparison and goodness-of-fit tests}\label{compare}

A consequence of Theorems \ref{E=R} and \ref{ED} is the construction of a statistical test of the hypothesis ${\cal H}_0: F=G$ against ${\cal H}_1: F\ne G$, based on two samples that may arise from correlated experiments. Let us choose the $b$-Wasserstein distance with $1<b<2$. The distributions $F$ and $G$ are supposed to be $C^2$ on $\mathbb R$ or $\mathbb R^+$ and satisfy $(CFG_{ED})$ and $(FG)$. By Theorem \ref{E=R}, under ${\cal H}_0$ the statistic $n^{b/2} W_{c}(\mathbb{F}_{n},\mathbb{G}_{n})$ converges to a positive finite random variable while by Theorem \ref{ED}, under ${\cal H}_1$ it converges almost surely to $+\infty$ at the rate $n^{b/2}W_{c}(F,G)$. Mathematically this test is effectively valid when the set $D=\{F^{-1}\ne G^{-1}\}$ is a finite union of non empty intervals, but we think that its validity could be extended to the more general case where $D$ is of positive Lebesgue measure in $(0,1)$. The use of $W_2^2$, with a rate $n$ is more restrictive since it needs very light tails. Nevertheless if sub-Gaussian tails can be asserted, by Theorem \ref{W2} the previous test works with $b=2$, which actually is a new test.

In each case the rather minimal $(CFG)$ type conditions have to be checked. They are close to be necessary in the proofs to overcome the difficulty of controlling how close the empirical tails of $\mathbb{F}_{n}$ and $\mathbb{G}_{n}$ must be under ${\cal H}_0$, and how far $\vert \mathbb{F}_{n}-\mathbb{G}_{n} \vert $ can deviate from $\vert {F}-{G} \vert $ in tails under ${\cal H}_1$. Interestingly the choice of $\rho(x)$ may be with a locally polynomial shape as $x \rightarrow 0$ and a different shape as $x \rightarrow +\infty$ possibly linear, polynomial or exponential. This flexibility allows to test the tail or the mid-quantiles with more or less accuracy.

In the same vein, concerning the distribution functions, Corollary \ref{pinball} yields
\begin{align*}
& \sqrt{n} \left( \int_{-\infty}^{+\infty} \left\vert \mathbb{F}_{n}(t)-\mathbb{G}_{n}(t) \right\vert dt - \int_{F^{-1}(D)} \left\vert F(t)-G(t) \right\vert dt \right) \\
& \rightarrow_{weak}\int_{D}\mathbb{B}(u)du+\int_{E}\left\vert \mathbb{B}%
(u)\right\vert du
\end{align*}
which seems not to have been already obtained. This provides weak limits for the power of the test under alternatives to ${\cal H}_0: F=G$ of the kind ${\cal H}_1: F=G_1$ where $G_1^{-1}$ only differs from $G^{-1}$ on an interval $D$, for instance with a slightly different right hand tail only. The test statistic $\sqrt{n} \int_{-\infty}^{+\infty} \left\vert \mathbb{F}_{n}(t)-\mathbb{G}_{n}(t) \right\vert dt$ has an almost sure first order rate of escape $\sqrt{n} \int_{G_1^{-1}(D)} \left\vert G_1(t)-G(t) \right\vert dt$.

As a by-product of the results of Section 3 one can similarly build goodness-of-fit tests ${\cal H}_0:F=F_0$ against ${\cal H}_1: F\ne F_0$ by using one sample under $F$ or by using an additional sample distributed as $F_0$. Notice that the test associated to $b=2$ was a consequence of \cite{Delbarrio05}.

\subsection{An application}\label{appli}

The motivation of our initial work was intimately related to the field of computer experiments. Many computer codes produce as output values of a function computed on so many points that it can be considered as a functional output. The case we are interested in is when this function is the $c.d.f.$ of a real $r.v.$ It turns out that Wasserstein distances are now commonly used to analyze such outputs. In view of defining new features for random $c.d.f.$ such as median or quantiles, more general Wasserstein costs may be used as contrasts to compute these features by solving an optimization problem -- see \cite{FK15}. Nevertheless computer codes only provide samples of the underlying distributions. Whence the importance of an efficient estimation of distances between $c.d.f.$ and goodness-of-fit tests through random samples.

As an illustration, let us conclude with a notion of quantile for a $r.v.$ taking values in the set of continuous $c.d.f.$'s. A useful new result of this article is the first part of Corollary \ref{pinball} which is strongly related to the preprint \cite{FK15}. Let $0<\alpha<1$. In \cite{FK15} the $\alpha-$quantile $F_\alpha$ of a random continuous $c.d.f.$ $\mathbb F$ is defined to be
$$F_\alpha=\argmin_{\theta\in \cal F} \mathbb E \ W_{c}(\mathbb F, \theta),$$
where $c(x,y)=(x-y) (\alpha-{\bf 1}_{x-y<0})$ is the non-symmetric contrast for  the $\alpha$-quantile of a real $r.v.$ and  $\cal F$ is the set of continuous $c.d.f.$ As previously mentioned, in practice a realization  $\mathbb F(\omega)$ of $\mathbb F$ is known through a $n$-sample of the distribution $\mathbb F(\omega)$. Hence we may assume that a $N$-sample  $\mathbb F_n^1$,\ldots,$\mathbb F_n^N$ is available, where each $\mathbb F_n^i$ is a $n$-empirical $c.d.f.$ of $\mathbb F^i$ and $\mathbb F^1$,\ldots,$\mathbb F^N$ are $i.i.d.$ according to $\mathbb F$. Define
 $$F_{N,n,\alpha}=\argmin_{\theta\in {\cal F}_n}  \frac{1}{N}\sum_{i=1}^NW_{c}(\mathbb F_n^i, \theta),$$
where ${\cal F}_n$ is the set of $c.d.f.$ with at most $n$ different values. Then one could use Corollary \ref{pinball} to prove that $F_{N,n,\alpha}$ is a consistent estimator of $F_\alpha$ when $N$ and $n$ tend to $+\infty$, and determine the rate of convergence. 

\section{Proofs}\label{sec:proof}

In the forthcoming proofs the high order quantiles are shown to have a secondary order impact compared to the mid-order quantiles that impose the rate as well as the limiting distribution under our sufficient conditions ensuring that the tails are not too heavy. For sake of simplicity we only work on the right hand tail, with quantiles of order $u\in\left(  \underline{u},1\right)$ for an arbitrary small $\underline{u}>0$. The counterpart for the left hand tail is immediate by using the same arguments.

To help the reader the variable of frequently used deterministic functions defined on $\mathbb{R}_{+}$ like $\rho_{\pm}$, $\rho_{\pm}^{-1}$, $l_{\pm}$, $l_{\pm}^{-1}$ or $L_{\pm}$ is denoted $x$ when considered as $x\rightarrow0$ and $y$ when considered as $y\rightarrow+\infty$. In the subsequent proofs the constant $K>0$ may change at each appearance.

In steps numbered 0 we remind active hypotheses while introducing local notation. The non standard Steps 1, 2 and 3 of the four proofs -- including the one in \cite{BFK17} -- are designed to address the non trivial problem of controlling the high order and extreme order quantiles under an explicit and almost minimal assumption on tails, namely $(CFG_{E})$, $(CFG_{D})$ or $(CFG_{ED})$. The secondary order terms in these conditions could be balanced slightly more sharply but at the price of adding technicalities to connect Steps 1 and 2. Finally we point out that the convergence at Steps 3 is weaker than in probability, due to the coupling approach.

\subsection{The case $F=G$}\label{proofF=G}
We prove Theorem \ref{E=R}.
\smallskip

\noindent\textbf{Step 0.} In this section $F=G$ and hence $E=\mathbb{R}$. For short, the key functions common to $X,Y$ are denoted $F^{-1}$, $\psi$, $H$ and $h$. Let assume $(FG)$, $(C)$ and $(CFG_{E})$ with $1\leqslant b_{\pm}<2$ in $(C2)$. Hence $\rho(x)=\max(\rho_{+}(x),\rho_{-}(x))\geqslant\rho_{c}(x)$ and $\rho_{\pm}(x)$ are positive convex increasing functions defined on $\mathbb{R}_{+}^{\ast}$ with $\rho_{\pm}(0)=0$. For $0\leqslant x\leqslant x_{0}$ we have $\rho_{\pm}\left(x\right)  =x^{b_{\pm}}L_{\pm}(x)$ and, whenever $b_{\pm}=1$ it is also assumed through (\ref{Lpi}) that $\lim_{x\rightarrow 0} L_{\pm}(x)=L_{\pm}(0)<+\infty$. Recall that $b=\min(b_{+},b_{-})$ and, for $0\leqslant x\leqslant x_{0}$, $\rho\left(  x\right)  =\max(\rho_{+}(x),\rho_{-}(x))=x^{b}L(x)$ where $L(x)$ is defined at (\ref{L}) and is slowly varying as $x\rightarrow 0$. We then have
\[
v_{n}=\frac{1}{\rho\left(  1/\sqrt{n}\right)  },\quad\lim_{n\rightarrow
+\infty}\frac{\sqrt{n}}{v_{n}}=1_{\{b=1\}}L(0)\text{.}%
\]
Since $L\in RV(0,0)$ we have, by the Karamata representation theorem,
\begin{equation}
L(x)=\exp\left(  \eta(x)+\int_{B}^{1/x}\frac{s(y)}{y}dy\right)  ,\quad 0<x\leqslant x_{0},\label{varlente}%
\end{equation}
with $B>0$, $\eta(x)$ and $s(y)$ are bounded measurable functions such that
\[
\lim_{x\rightarrow0}\eta(x)=\eta_{\infty}\in\mathbb{R},\quad\lim
_{y\rightarrow+\infty}s(y)=0.
\]
We can then define
\begin{equation}
\eta_{0}=\sup_{0<x\leqslant x_{0}}\left\vert \eta(x)\right\vert \in\mathbb{R}_+ , \quad c_0=e^{2\eta_{0}}\geqslant 1.\label{eta0}
\end{equation}
For $y$ large it holds $\rho_{\pm}\left(  y\right)  =\exp(l_{\pm}(y))$ where the functions $l_{\pm}(y)$ are not asked to be in$\ {RV}(+\infty,\gamma_{\pm})$ in this proof, but (\ref{L'}) does matter. However in practice if $(C3)$ would not hold then $(CFG_{E})$ would be more difficult to translate in terms of admissible $F$. Hence, for some $y_{0} > x_{0}$,
\[
\rho\left(  y\right)  =\exp(l(y)),\quad l(y)=\max(l_{+}(y),l_{-}(y)),\quad
y \geqslant y_{0}.
\]
Since $\rho_{\pm}$ and $\rho$ are convex, by (\ref{L'}) there exists $d_{\pm}\geqslant1$, $d=\min(d_{-},d_{+})$ and $d_{0,\pm}$, $d_{0}$ such that%
\begin{equation}
l_{\pm}(y)\geqslant d_{\pm}\log y+d_{0,\pm},\quad l(y)\geqslant d\log y+d_{0},\quad y \geqslant y_{0}.\label{dlogy}%
\end{equation}
By $(CFG_{E})$, the joint influence of $l$, $L$ and $b$ on the allowed tails $F^{-1}$ is expressed at (\ref{lrondpsy}) if $b>1$ and (\ref{lrondpsyb1}) if $b=1$.

We decompose the integral $W_{c}(\mathbb{F}_{n},\mathbb{G}_{n})$ as follows, with the three remainder terms implicitly treated in a similar way for left hand tails. We will specify later two positive sequences $i_{n}$ and $j_{n}$ such that $n>j_{n}>i_{n}\rightarrow+\infty$.\ The proof consists in four steps, each dealing with one of the four terms
\begin{equation}
W_{c}(\mathbb{F}_{n},\mathbb{G}_{n})=I_{\mathcal{I}_{n}}+I_{\mathcal{J}%
_{n}}+I_{\mathcal{K}_{n}}+I_{\mathcal{L}},\quad I_{A}=\int_{A}\rho_{c}\left(
\mathbb{F}_{n}^{-1}(u)-\mathbb{G}_{n}^{-1}(u)\right)  du,\label{WcE}%
\end{equation}
where $\mathcal{I}_{n}=\left(  1-i_{n}/n,1\right]  $, $\mathcal{J}_{n}=\left(1-j_{n}/n,1-i_{n}/n\right]  $, $\mathcal{K}_{n}=\left(  \overline{u}, 1-j_{n}/n\right]  $, 
$\mathcal{L}=\left[  \underline{u},\overline{u}\right]$ and $0<\underline{u}<1/2<\overline{u}<1$. In order to accurately choose $i_{n}$ and $j_{n}$ one has to take into account two difficulties. First, the rate $1/v_{n}$\ is faster than $1/\sqrt{n}$ so that $\mathcal{I}_{n}\cup\mathcal{J}_{n}$ should be sufficiently small. Second, the empirical extreme quantile difference $\mathbb{F}_{n}^{-1}(u)-\mathbb{G}_{n}^{-1}(u)$ may be either very large or\ very small as $u\rightarrow1$, thus the cost function $\rho_{c}(\mathbb{F}_{n}^{-1}(u)-\mathbb{G}_{n}^{-1}(u))$ is evaluated at $0$ on some random subsets of $\mathcal{I}_{n}\cup\mathcal{J}_{n}$ and at $+\infty$ on some others. The later problem is the most difficult to address.\textit{\smallskip}

\noindent\textbf{Step 1.} Let $K_{n}$ be a positive sequence such that $K_{n}\rightarrow+\infty$ and define
\begin{equation}
i_{n}=\frac{n}{v_{n}K_{n}\rho(\psi^{-1}(\log n+K_{n}))}. \label{in}%
\end{equation}
Notice that $(FG1)$ and (\ref{dlogy}) imply that $\rho(\psi^{-1}(\log n+K_{n}))\rightarrow+\infty$ and $i_{n}=o\left(  n^{1-b/2}L(1/\sqrt{n})/K_{n}\right)$ as $n\rightarrow+\infty$, so that $i_{n}/ \sqrt{n}\rightarrow0$ even when $b=1$, thanks to (\ref{Lpi}). The following lemma ensures that $i_{n}{/\log\log n\rightarrow+\infty}$. Observe also that $\psi^{-1}(\log n+K_{n})=F^{-1}(1-1/ne^{K_{n}})$ is an extreme quantile just beyond the expected order $F^{-1}(1-1/n)$ for $X_{(n)}$ and $Y_{(n)}$, which is the key to Lemma \ref{Lem_in}. Let $[y]$ denote the integer part of $y$. Consider the $r.v.$
\[
I_{\mathcal{I}_{n}}\leqslant\int_{\mathcal{I}_{n}}\rho\left(  \mathbb{F}_{n}^{-1}(u)-\mathbb{G}_{n}^{-1}(u)\right)  du=\frac{1}{n}{\sum\limits_{i=n-[i_{n}]}^{n}}\rho\left(X_{(i)}-Y_{(i)}\right).
\]

\begin{lemma}
\label{Lem_in}Assume $(FG1)$, $(C)$ and $(CFG_{E})$. There exists $K_{n}$ such
that
\[
K_{n}\rightarrow+\infty,\quad\underset{n\rightarrow+\infty}{\lim}\frac{K_{n}%
}{{\log\log n}}{=0},\quad\underset{n\rightarrow+\infty}{\lim\inf}\frac{\log
i_{n}}{{\log\log n}}{\geqslant\theta _2>0}%
\]
and%
\[
\lim_{n\rightarrow+\infty}v_{n}I_{\mathcal{I}_{n}}=0\;\text{in probability.}%
\]
\end{lemma}

\noindent\textbf{Proof.} \textbf{(i)} Let $K_{n}\rightarrow+\infty$, $K_{n}/\log\log n\rightarrow0$ be as slow as needed later. By $(FG1)$ we have $F^{-1}\left(  1-1/ne^{K_{n}}\right)  \rightarrow+\infty$ as $n\rightarrow+\infty$, yet arbitrarily slowly. Thus, by (\ref{lrondpsy}) and (\ref{in}) we have, for any $\theta^{\prime\prime}>1-b/2$, any $\theta^{\prime}<\theta _2$ and all $n$ large,
\begin{align*}
i_{n}  &  =\frac{n^{1-b/2}L(1/\sqrt{n})}{K_{n}\rho(\psi^{-1}(\log n+K_{n}))}\\
&  \geqslant\frac{1}{K_{n}}\frac{L(1/\sqrt{n})}{L\left(  1/\sqrt{ne^{K_{n}}%
}\right)  }\exp\left(  -\left(  1-\frac{b}{2}\right)  K_{n}+2\log\psi
^{-1}(\log n+K_{n})+\theta _2\log(\log n+K_{n})\right) \\
&  \geqslant\frac{L(1/\sqrt{n})}{L\left(  1/\sqrt{ne^{K_{n}}}\right)  }%
\frac{1}{e^{\theta^{\prime\prime}K_{n}}}\left(  F^{-1}\left(  1-\frac
{1}{ne^{K_{n}}}\right)  \right)  ^{2}(\log n+K_{n})^{\theta _2}\\
&  \geqslant\frac{L(1/\sqrt{n})}{L\left(  1/\sqrt{ne^{K_{n}}}\right)  }(\log
n)^{\theta^{\prime}}.
\end{align*}
Applying (\ref{varlente}) and $K_{n}\rightarrow+\infty$ we get%
\[
\frac{L\left(  1/\sqrt{ne^{K_{n}}}\right)  }{L(1/\sqrt{n})}=\exp\left(
\eta\left(  \frac{1}{\sqrt{ne^{K_{n}}}}\right)  -\eta\left(  \frac{1}{\sqrt
{n}}\right)  +\int_{\sqrt{n}}^{\sqrt{ne^{K_{n}}}}\frac{s(y)}{y}dy\right) .
\]
Since $e^{K_{n}}<\log n$ we can furthermore choose $K_{n}$ such that%
\[
K_{n}<\frac{1}{s_{n}},\quad s_{n}=\sup_{\sqrt{n}\leqslant y\leqslant
\sqrt{n\log n}}s(y),
\]
where $s_{n}\rightarrow0$ as $n\rightarrow+\infty$. The slower is $L$ the
faster is $1/s_{n}$ hence the resulting requirement is sometimes only the
initial $K_{n}/\log\log n\rightarrow0$. We readily obtain, by (\ref{eta0}),
\[
\underset{n\rightarrow+\infty}{\lim\sup}\ \frac{L\left(  1/\sqrt{ne^{K_{n}}%
}\right)  }{L(1/\sqrt{n})}\leqslant\underset{n\rightarrow+\infty}{\lim\sup
}\ \exp\left(  2\eta_{0}+s_{n}\frac{K_{n}}{2}\right)  <+\infty.
\]
The claimed deterministic $\lim\inf$ is proved by letting $\theta^{\prime} \rightarrow \theta _2$. Notice that $(CFG_{E})$ was crucially required.\smallskip

\noindent\textbf{(ii)} Concerning the stochastic integral $I_{\mathcal{I}_{n}%
}$ the choice of $i_{n}$ in (\ref{in}) is minimal to guaranty the rate $v_{n}%
$ and $(CFG_{E})$ is not required. Recall that $F$ has support $\mathbb{R}$. Fix $\varepsilon>0$ and consider the events
\[
A_{n}=\left\{  v_{n}I_{\mathcal{I}_{n}}\geqslant4\varepsilon\right\}  ,\quad
B_{n}=\left\{  X_{(n-[i_{n}])}>0\cap Y_{(n-[i_{n}])}>0\right\}  .
\]
We have $\mathbb{P}\left(  A_{n}\right)  \leqslant\mathbb{P}\left(  A_{n}\cap
B_{n}\right)  +\mathbb{P}\left(  B_{n}^{c}\right)  $ and $\mathbb{P}\left(
B_{n}^{c}\right)  \rightarrow0$ as $n \rightarrow +\infty$. On $B_{n}$ it holds%
\[
v_{n}I_{\mathcal{I}_{n}}\leqslant\frac{v_{n}}{n}{\sum\limits_{i=n-[i_{n}]}%
^{n}}\left(  \rho_{+}\left(  X_{(i)}\right)  +\rho_{-}\left(  Y_{(i)}\right)
\right)  \leqslant\frac{v_{n}}{n}(i_{n}+1)\left(  \rho_{+}\left(
X_{(n)}\right)  +\rho_{-}\left(  Y_{(n)}\right)  \right)
\]
hence $\mathbb{P}\left(  A_{n}\cap B_{n}\right)  \leqslant\mathbb{P}\left(
C_{n,X}\right)  +\mathbb{P}\left(  C_{n,Y}\right)  $ where%
\[
C_{n,X}=\left\{  \rho_{+}\left(  X_{(n)}\right)  \geqslant\varepsilon\frac
{n}{v_{n}i_{n}}\right\}  ,\quad C_{n,Y}=\left\{  \rho_{-}\left(
Y_{(n)}\right)  \geqslant\varepsilon\frac{n}{v_{n}i_{n}}\right\}  .
\]
In order to evaluate $\mathbb{P}\left(  C_{n,X}\right)  =1-\left(
1-\mathbb{P}\left(  \rho_{+}(X)>\varepsilon n/v_{n}i_{n}\right)  \right)
^{n}$ we combine $\rho_{+}^{-1}(x)=l_{+}^{-1}(\log x)$, $l^{-1}\leqslant
l_{+}^{-1}$ and $\psi_{X}=\psi$ with (\ref{in}) to obtain, for $n$ large
enough to have $K_{n}>1/\varepsilon$,%
\begin{align*}
&  \mathbb{P}\left(  \rho_{+}(X)>\varepsilon K_{n}\rho(\psi^{-1}(\log
n+K_{n}))\right) \\
&  \leqslant\exp\left(  -\psi\circ l^{-1}\left(  \log\varepsilon+\log
K_{n}+l(\psi^{-1}(\log n+K_{n}))\right)  \right) \\
&  \leqslant\frac{1}{ne^{K_{n}}}.
\end{align*}
Therefore $\mathbb{P}\left(  C_{n,X}\right)  \leqslant1-\exp\left(
-\exp(-K_{n})\right)  \sim\exp(-K_{n})\rightarrow0$ as $n\rightarrow+\infty$,
and similarly $\mathbb{P}\left(  C_{n,Y}\right)  \rightarrow0$. This implies
that $v_{n}I_{\mathcal{I}_{n}}\rightarrow0$ in probability.$\quad
\square\smallskip$

\noindent\textbf{Step 2.}\ Write $\beta_{n}(u)=\beta_{n}^{X}(u)-\beta_{n}%
^{Y}(u)$ with%
\begin{equation}
\beta_{n}^{X}(u)=\sqrt{n}(\mathbb{F}_{n}^{-1}(u)-F^{-1}(u)),\quad\beta_{n}%
^{Y}(u)=\sqrt{n}(\mathbb{G}_{n}^{-1}(u)-F^{-1}(u)), \label{quantile}%
\end{equation}
thus $I_{A}=\int_{A}\rho_{c}\left(  \beta_{n}(u)/\sqrt{n}\right)  du$ in
(\ref{WcE}). Let $\Delta_{n}=\mathcal{J}_{n}\cup\mathcal{K}_{n}\cup
\mathcal{L}=\left[  \underline{u},1-i_{n}/n\right]  $. The next lemma shows
that in the integral $I_{\Delta_{n}}$ the cost function $\rho$ is evaluated
near $0$ provided that $n$ is large.

\begin{lemma}
\label{Lem_tube}Assume $(FG)$ and $(CFG_{E})$. For any $0<\xi<1/2-b/4$ it
holds%
\[
\underset{n\rightarrow+\infty}{\lim}(\log n)^{\xi}\sup_{u\in\Delta_{n}}%
\frac{\left\vert \beta_{n}(u)\right\vert }{\sqrt{n}}={0}\quad\text{a.s.}%
\]

\end{lemma}

\noindent\textbf{Proof.} \textbf{(i)} Assuming $(FG1)$, $(FG2)$ and since
$i_{n}/\log\log n\rightarrow+\infty$ by Lemma \ref{Lem_in} we can apply the
classical hungarian results to $\left\vert \beta_{n}(u)\right\vert
\leqslant\left\vert \beta_{n}^{X}(u)\right\vert +\left\vert \beta_{n}%
^{Y}(u)\right\vert $ exactly as for Lemma 23 in \cite{BFK17} to
get%
\begin{equation}
\underset{n\rightarrow+\infty}{\lim\sup}\sup_{u\in\Delta_{n}}\frac
{h(u)\left\vert \beta_{n}(u)\right\vert }{\sqrt{\left(  1-u\right)  \log\log
n}}\leqslant8\quad\text{a.s.} \label{hongr}%
\end{equation}
Next observe that $(FG3)$ implies, for some $0<M<+\infty$ and $u\in\Delta_{n}%
$,%
\begin{align}
\frac{1}{M}\frac{\sqrt{1-u}}{h(u)}\sqrt{\frac{\log\log n}{n}}  &
\leqslant\frac{F^{-1}(u)}{\sqrt{1-u}}\sqrt{\frac{\log\log n}{n}}\nonumber\\
&  \leqslant\varepsilon_{n}=\frac{F^{-1}(1-i_{n}/n)}{\sqrt{i_{n}}}\sqrt
{\log\log n}. \label{epsiln}%
\end{align}
\textbf{(ii)} Remind that $e^{-K_{n}}<1<i_{n}$ for all $n$ large, and
$F^{-1}\left(  1-i_{n}/n\right)  \rightarrow+\infty$ as $n\rightarrow+\infty$
with no obvious control on the rate. By (\ref{in}) and the consequence
(\ref{lrondpsy}) of $(CFG_{E})$ we have already seen in the proof of Lemma
\ref{Lem_in} that if $\theta^{\prime\prime}<1-b/2$ and $\theta^{\prime\prime
}<\theta^{\prime}<\theta _2$ then it holds, for all $n$ large enough,%
\begin{align*}
i_{n}  &  \geqslant\frac{1}{e^{\theta^{\prime\prime}K_{n}}}\left(
F^{-1}\left(  1-\frac{1}{ne^{K_{n}}}\right)  \right)  ^{2}(\log n+K_{n}%
)^{\theta _2}\\
&  \geqslant\left(  F^{-1}\left(  1-\frac{1}{n}\right)  \right)  ^{2}(\log
n)^{\theta^{\prime}}\\
&  \geqslant\left(  F^{-1}\left(  1-\frac{i_{n}}{n}\right)  \right)  ^{2}%
(\log\log n)(\log n)^{\theta^{\prime\prime}}%
\end{align*}
hence for any $0<\xi<\theta^{\prime\prime}/2$ it holds $\lim_{n\rightarrow
+\infty}(\log n)^{\xi}\varepsilon_{n}=0.$ The conclusion follows, by
(\ref{hongr}) and (\ref{epsiln}).$\quad\square\smallskip$

Let $j_{n}=n^{\beta}$ with $1/2<\beta<1$, so that $i_n<\sqrt{n}<j_n$ for all $n$ large. Remind $\varepsilon_{n}$ from (\ref{epsiln}). Let introduce%
\begin{equation}
\varepsilon_{n}(u)=9\frac{\sqrt{1-u}}{h(u)}\sqrt{\frac{\log\log n}{n}%
}\leqslant9\varepsilon_{n},\quad u\in\mathcal{J}_{n}. \label{epsilnu}%
\end{equation}

\begin{lemma}
\label{Lem_E_step2}Assume $(FG)$, $(C)$ and $(CFG_{E})$. Then we have
\[
\lim_{n\rightarrow+\infty}v_{n}I_{\mathcal{J}_{n}}=0\quad a.s.
\]

\end{lemma}

\noindent\textbf{Proof.} \textbf{(i)} By Lemma \ref{Lem_tube}, for all $n$
large enough and any $u\in\mathcal{J}_{n}$ it holds%
\[
\frac{1}{n}\leqslant\sup_{u\in\mathcal{J}_{n}}\sqrt{\frac{1-u}{n}}%
\leqslant\varepsilon_{n}(u)\leqslant\frac{1}{(\log n)^{\xi}}.
\]
Consider $L$ defined in (\ref{L}). Using (\ref{varlente}) and (\ref{eta0}) we get%
\begin{equation}
L_{n}=\sup_{u\in\mathcal{J}_{n}}\frac{L(\varepsilon_{n}(u))
}{L(1/\sqrt{n}) }\leqslant\exp\left(  2\eta_{0}+%
{\displaystyle\int\nolimits_{(\log n)^{\xi}}^{n}}
\frac{\left\vert s(y)\right\vert }{y}dy\right)  \label{Ln}%
\end{equation}
hence%
\begin{equation}
\lim_{n\rightarrow+\infty}\frac{\log L_{n}}{\log n}\leqslant\lim
_{n\rightarrow+\infty}\frac{1}{\log n}\left(  2\eta_{0}+\log n\sup_{(\log
n)^{\xi}\leqslant y\leqslant n}\left\vert s(y)\right\vert \right)  =0.
\label{logLn}%
\end{equation}
\textbf{(ii)} Remind that $\rho_{\pm}$ are increasing. By Lemma \ref{Lem_tube}
and $(C2)$ we almost surely have, for all $n$ large,%
\[
I_{\mathcal{J}_{n}}\leqslant\int_{\mathcal{J}_{n}\cap\left\{  \beta
_{n}\geqslant0\right\}}\rho_{+}(\varepsilon_{n}(u))
du+\int_{\mathcal{J}_{n}\cap\left\{  \beta_{n}<0\right\}  }\rho_{-}(
\varepsilon_{n}(u))du
\]
where, by (\ref{hongr}), (\ref{epsiln}) and (\ref{epsilnu}), $\sup
_{u\in\mathcal{J}_{n}}\varepsilon_{n}(u)\leqslant9\varepsilon_{n}\rightarrow
0$. Hence, recalling (\ref{L}) we are reduced to study the bounding
deterministic integral%
\[
I_{\mathcal{J}_{n}}\leqslant\int_{\mathcal{J}_{n}}\rho(  \varepsilon
_{n}(u))  du=\int_{\mathcal{J}_{n}}\left(  \varepsilon_{n}(u)\right)
^{b}L(\varepsilon_{n}(u)) du.
\]
By (\ref{L}), $L_{n}$ from (\ref{Ln}) and $(FG3)$ we further have%
\begin{equation}
v_{n}I_{\mathcal{J}_{n}}\leqslant L_{n}(\log\log n)^{b/2}\int_{\mathcal{J}%
_{n}}\left(  \frac{F^{-1}(u)}{\sqrt{1-u}}\right)  ^{b}du. \label{vnIJn}%
\end{equation}
We next show that $L_{n}(\log\log n)^{b/2}$ is a secondary order factor
compared to the integral in (\ref{vnIJn}), whatever the choice of $1/2<\beta<1$ defining $j_n$ in $\mathcal{J}_{n}$.\smallskip

\noindent\textbf{(iii)} The fact that $l(y)\geqslant\log y$ as $y\rightarrow
+\infty$ combined to $(CFG_{E})$ shows that for all $u$ large enough, we have%
\begin{align*}
F^{-1}(u)  &  =\psi^{-1}\left(  \log\frac{1}{1-u}\right) \\
&  \leqslant\exp\left(  l\circ\psi^{-1}\left(  \log\frac{1}{1-u}\right)
\right) \\
&  \leqslant\exp\left(  \left(  1-\frac{b}{2}\right)  \log\frac{1}{1-u}+\log
L(  \sqrt{1-u})  -2\log F^{-1}(u)-\theta _2\log\log\frac{1}%
{1-u}\right).
\end{align*}
Therefore we get%
\[
F^{-1}(u)^{3}\leqslant\left(  \frac{1}{1-u}\right)  ^{1-b/2}\frac{L(\sqrt
{1-u})}{(\log(1/(1-u)))^{\theta_2}}%
\]
and%
\begin{align}
\int_{\mathcal{J}_{n}}\left(  \frac{F^{-1}(u)}{\sqrt{1-u}}\right)  ^{b}du  &
\leqslant\int_{\mathcal{J}_{n}}\left(  \frac{1}{1-u}\right)  ^{(1-b/2)b/3+b/2}%
\frac{L(\sqrt{1-u})^{b/3}}{(\log(1/(1-u)))^{\theta_2 b/3}}du\nonumber\\
&  \leqslant\int_{\mathcal{J}_{n}}\left(  \frac{1}{1-u}\right)  ^{b(5-b)/6}%
\frac{L(\sqrt{1-u})^{b/3}}{(\log(1/(1-u)))^{\theta_2 b/3}}du. \label{Lb3}%
\end{align}
Since $1\leqslant b<2$ we can find $\gamma$ such that $0<b(5-b)/6<\gamma<1$.
The second factor in the integral (\ref{Lb3}) is slowly varying in $1-u$ as
$u\rightarrow1$ thus the whole integral is ultimately bounded from above by%
\begin{equation}
(1-\gamma)\int_{\mathcal{J}_{n}}\left(  \frac{1}{1-u}\right)  ^{\gamma
}du=\left[  -(1-u)^{1-\gamma}\right]  _{1-j_{n}/n}^{1-i_{n}/n}\leqslant
\frac{1}{n^{(1-\gamma)(1-\beta)}}. \label{integamma}%
\end{equation}
We deduce from (\ref{logLn}), (\ref{vnIJn}), (\ref{Lb3}) and (\ref{integamma}) the convergence
\[
\lim_{n\rightarrow+\infty}v_{n}I_{\mathcal{J}_{n}}\leqslant\lim_{n\rightarrow
+\infty}\frac{L_{n}(\log\log n)^{b/2}}{n^{(1-\gamma)(1-\beta)}}=0\quad a.s.
\]
at a power rate.$\quad\square\smallskip$

\noindent\textbf{Step 3.}\ Compared to $\mathcal{J}_{n}$ the interval
$\mathcal{K}_{n}$ is so large that $v_{n}I_{\mathcal{K}_{n}}$ can no more
converge to zero in probability. Instead it is made small with high probability by choosing
$\overline{u}$ and $\beta$ properly, at Lemma \ref{Lem_E_step3}. Moreover, in
order to evaluate the integral of $\rho_{c}(\beta_{n}(u)/\sqrt{n})$ over
$\mathcal{K}_{n}$ accurately enough it is no more sufficient to bound the
process, therefore we approximate it at Lemma \ref{Lem_E_approxim} by a
Gaussian process which helps revealing the underlying deterministic integral
to compute. Lastly the fact that $\beta_{n}(u)$ itself may be very small or
very large along $\mathcal{K}_{n}$ makes a bit tedious the uniform control of
the slowly varying part $L(x)$ of $\rho(x)$.

Define $\Delta_{n}^{\prime}=\left(  j_{n}/n,1-j_{n}/n\right)  $. We
first recall the strong approximation of the joint quantile processes
\[
\mathcal{Q}_{n}(u)=\left(\beta_{n}^{X}(u),\beta_{n}^{Y}(u)\right), \quad u\in\Delta_{n}^{\prime},
\]
by the joint Gaussian processes%
\[
\mathcal{G}_{n}(u)= \left( \mathbb{B}_{n}^{X}(u) ,\mathbb{B}_{n}^{Y}(u)\right),\quad \mathbb{B}_{n}^{X}(u)=\frac{B_{n}^{X}(u)}{h_X(u)},\quad \mathbb{B}_{n}^{Y}(u)=\frac{B_{n}^{Y}(u)}{h_Y(u)},\quad u\in\Delta_{n}^{\prime},
\]
where $B_{n}^{X}(u)=\mathbb{H}_{n}(H_{X}(u))$, $B_{n}^{Y}(u)=\mathbb{H}_{n}(H^{Y}(u))$ and $\mathbb{H}_{n}$ is a $\mathbb{P}^{X,Y}$-Brownian bridge indexed by the halfplanes
\[
H_{X}(u)=\left\{  (x,y):x\leqslant F^{-1}(u)\right\}  ,\quad H^{Y}(u)=\left\{
(x,y):y\leqslant F^{-1}(u)\right\}.
\]
Therefore $B_{n}^{X}$ and $B_{n}^{Y}$ are two standard Brownian bridges with cross covariance given for $u,v \in (0,1)$ by
\begin{align*}
cov(B_{n}^{X}(u),B_{n}^{Y}(v)) & = \mathbb{P}^{X,Y}(H_{X}(u)\cap H^{Y}(v))-\mathbb{P}^{X,Y}(H_{X}(u))\mathbb{P}^{X,Y}(H^{Y}(v)) \\
& = \mathbb{P}\left(X\leqslant F^{-1}(u),Y\leqslant F^{-1}(v)\right) - uv\\
& = H(F^{-1}(u),F^{-1}(v)) - uv.
\end{align*}
Notice that $H(F^{-1}(u),F^{-1}(v))$ is the copula function of $(X,Y)$. From now and for the remainder of the proof
we work on the probability space of the following Lemma \ref{Lem_E_approxim}.
The weak convergence finally established on this space at Steps 4 and 5 remains valid
on any probability space.$\smallskip$

\begin{lemma}
\label{Lem_E_approxim}Assume $(FG)$. Then we can build on the same probability
space versions of $\left(  X_{n},Y_{n}\right)  _{n\geqslant1}$ and $\left(
\mathbb{H}_{n}\right)  _{n\geqslant1}$ such that $\mathcal{Q}_{n}(u)=$
$\mathcal{G}_{n}(u)+\mathcal{Z}_{n}(u)$ for all $n\geqslant1$ and
$u\in\Delta_{n}^{\prime}$ where $\mathcal{Z}_{n}(u)=(Z_{n}^{X}(u)/h_X(u),Z_{n}^{Y}(u)/h_Y(u))$ satisfies, for some $\upsilon\in\left(  0,1/22\right)  $,%
\[
\lim_{n\rightarrow+\infty}n^{\upsilon}\sup_{u\in\Delta_{n}^{\prime}}\left\vert
Z_{n}^{X}(u)\right\vert =\lim_{n\rightarrow+\infty}n^{\upsilon}\sup
_{u\in\Delta_{n}^{\prime}}\left\vert Z_{n}^{Y}(u)\right\vert =0\quad a.s.
\]

\end{lemma}

\noindent\textbf{Proof.} This follows from Theorem 28 in \cite{BFK17} with $F=G$.$\quad\square\smallskip$

The joint strong approximation of Lemma \ref{Lem_E_approxim} applied with $F=G$ and $h_X=h_Y=h$ combined to $(CFG_{E})$ provides a stochastic control of the deviations of $v_{n}I_{\mathcal{K}_{n}}$ that is
weaker than in probability but sufficient for the targeted weak convergence. Since it concerns the probability distribution of $I_{\mathcal{K}_{n}}$ the following lemma remains true on any probability space.$\smallskip$

\begin{lemma}
\label{Lem_E_step3}Assume $(FG)$, $(C)$ and $(CFG_{E})$. There exists
$\beta\in\left(  1/2,1\right)  $ such that for any choice of $\lambda>0$ and
$\varepsilon>0$ one can find $\overline{u}_0\in\left(  1/2,1\right)  $ and
$n_{0}>0$ such that, for all $\overline{u} \in [\overline{u}_0,1)$ and $n>n_{0}$,
\[
\mathbb{P}\left(  v_{n}I_{\mathcal{K}_{n}}>\lambda\right)  <\varepsilon
\text{.}%
\]

\end{lemma}

\noindent\textbf{Proof.} Fix $\lambda>0$ and $\varepsilon>0$ then consider, with $\beta_{n}$ as in (\ref{quantile}) the event
\[
\mathcal{C}_{n}^{\lambda}=\left\{  v_{n}\int_{\mathcal{K}_{n}}\rho_{c}\left(
\frac{\beta_{n}(u)}{\sqrt{n}}\right)  du>\lambda\right\}  .
\]
\textbf{(i)} For $0<\tau<\min(1,\lambda/2)$ define the random sets
\[
\mathcal{K}_{n}^{<\tau}=\left\{  u\in\mathcal{K}_{n}:\left\vert \beta
_{n}(u)\right\vert <\tau\right\}  ,\quad\mathcal{K}_{n}^{>\tau}=\mathcal{K}%
_{n}\backslash\mathcal{K}_{n}^{<\tau}.
\]
Recalling that the cost $\rho_{\pm}$ is convex, positive and such that
$\rho_{\pm}(0)=0$\ we have $\rho_{\pm}(\tau x)\leqslant\tau\rho_{\pm}(x)$ for
all $x\geqslant0$. It follows that%
\begin{align*}
v_{n}I_{\mathcal{K}_{n}^{<\tau}}  &  \leqslant v_{n}\int_{\mathcal{K}%
_{n}^{<\tau}\cap\left\{  \beta_{n}<0\right\}  }\tau\rho_{+}\left(  \frac
{1}{\sqrt{n}}\right)  du+v_{n}\int_{\mathcal{K}_{n}^{<\tau}\cap\left\{
\beta_{n}\geqslant0\right\}  }\tau\rho_{-}\left(  \frac{1}{\sqrt{n}}\right)
du\\
&  \leqslant\frac{\max(\rho_{-}\left(  1/\sqrt{n}\right)  ,\rho_{+}\left(
1/\sqrt{n}\right)  )}{\rho\left(  1/\sqrt{n}\right)  }\tau\int_{\mathcal{K}%
_{n}^{<\tau}}du\leqslant\tau.
\end{align*}
As a consequence,%
\[
\mathbb{P}\left(  \mathcal{C}_{n}^{\lambda}\right)  =\mathbb{P}\left(
v_{n}(I_{\mathcal{K}_{n}^{<\tau}}+I_{\mathcal{K}_{n}^{>\tau}})>\lambda\right)
\leqslant\mathbb{P}\left(  v_{n}I_{\mathcal{K}_{n}^{>\tau}}\geqslant
\lambda-\tau\right)  \leqslant\mathbb{P}\left(  v_{n}I_{\mathcal{K}_{n}%
^{>\tau}}\geqslant\frac{\lambda}{2}\right)  .
\]
\textbf{(ii)} For all $n\geqslant n_{0}$ and $n_{0}=n_{0}(\varepsilon,\xi)$
large enough we have $(\log n)^{\xi}<\sqrt{n}$ together with, by Lemma
\ref{Lem_tube} and since $\mathcal{K}_{n}^{>\tau}\subset \mathcal{K}_{n}\subset\Delta_{n}$,%
\[
\mathbb{P}\left(  \mathcal{D}_{n}\right)  >1-\frac{\varepsilon}{2}%
,\quad\mathcal{D}_{n}=\left\{  \sup_{u\in\mathcal{K}_{n}^{>\tau}}%
\frac{\left\vert \beta_{n}(u)\right\vert }{\sqrt{n}}\leqslant\frac{1}{(\log
n)^{\xi}}\right\}  .
\]
Assume now that $n\geqslant n_{0}$. On the event $\mathcal{D}_{n}$, for any
$u\in\mathcal{K}_{n}^{>\tau}$ we have%
\[
\frac{\tau}{\sqrt{n}}\leqslant\min\left(  \frac{1}{\sqrt{n}},\frac{\left\vert
\beta_{n}(u)\right\vert }{\sqrt{n}}\right)  \leqslant\frac{1}{(\log n)^{\xi}}%
\]
which by (\ref{L}), (\ref{varlente}) and (\ref{eta0}) yields
\begin{align*}
\frac{L\left(  \left\vert \beta_{n}(u)\right\vert /\sqrt{n}\right)  }{L\left(
1/\sqrt{n}\right)  }  &  \leqslant\exp\left(  2\eta_{0}+%
{\displaystyle\int\nolimits_{\min(\sqrt{n},\sqrt{n}/\left\vert \beta
_{n}(u)\right\vert )}^{\sqrt{n}/\tau}}
\frac{\left\vert s(y)\right\vert }{y}dy\right) \\
&  \leqslant c_{0}\exp\left(  s_{n}%
{\displaystyle\int\nolimits_{\min(\sqrt{n},\sqrt{n}/\left\vert \beta
_{n}(u)\right\vert )}^{\sqrt{n}/\tau}}
\frac{1}{y}dy\right) \\
&  =c_{0}\exp\left(  s_{n}\left(  \max\left(  0,\log(\left\vert \beta
_{n}(u)\right\vert )\right)  -\log\tau\right)  \right) \\
&  \leqslant c_{0}\left\vert \frac{\beta_{n}(u)}{\tau}\right\vert ^{q_{n}(u)}%
\end{align*}
where the sequence $s_{n}$ and the stochastic process $q_{n}(u)$ are defined by
\begin{equation}
s_{n}=\sup_{(\log n)^{\xi}\leqslant y\leqslant\sqrt{n}/\tau}\left\vert
s(y)\right\vert ,\quad q_{n}(u)=s_{n}1_{\left\{  \left\vert \beta
_{n}(u)\right\vert >1\right\}  }. \label{sn}%
\end{equation}
Since $s(y)\rightarrow0$ as $y\rightarrow+\infty$ we have%
\begin{equation}
\lim_{n\rightarrow+\infty}\sup_{u\in\mathcal{K}_{n}^{>\tau}}q_{n}%
(u)\leqslant\lim_{n\rightarrow+\infty}s_{n}=0 \label{qnu0}%
\end{equation}
and this uniform convergence of $q_{n}$ is certain, not almost sure. In other
words, the uncertainty in the following inequality only comes from
$\mathbb{P}\left(  \mathcal{D}_{n}\right)  $. We have shown that for all $n$
large, on the event $\mathcal{D}_{n}$, it holds%
\begin{equation}
v_{n}I_{\mathcal{K}_{n}^{>\tau}}\leqslant v_{n}\int_{\mathcal{K}_{n}^{>\tau}%
}\rho\left(  \frac{\left\vert \beta_{n}(u)\right\vert }{\sqrt{n}}\right)
du\leqslant\frac{c_{0}}{\tau^{s_{n}}}\int_{\mathcal{K}_{n}^{>\tau}}\left\vert
\beta_{n}(u)\right\vert ^{b+q_{n}(u)}du \label{IKnbqn}%
\end{equation}
where $\tau^{s_{n}}\rightarrow1$ as $n\rightarrow+\infty$. We are ready to bound
$\mathbb{P}\left(  \mathcal{D}_{n}\cap\left\{  v_{n}I_{\mathcal{K}_{n}^{>\tau
}}\geqslant\lambda/2\right\}  \right)  $.$\smallskip$

\noindent\textbf{(iii)} On the probability space of Lemma \ref{Lem_E_approxim}
we have%
\[
\left\vert \beta_{n}(u)\right\vert \leqslant\frac{\left\vert B_{n}%
^{X}(u)\right\vert }{h(u)}+\frac{\left\vert B_{n}^{Y}(u)\right\vert }%
{h(u)}+\frac{\left\vert Z_{n}^{X}(u)\right\vert }{h(u)}+\frac{\left\vert
Z_{n}^{Y}(u)\right\vert }{h(u)}.
\]
If $\alpha\geqslant1$ then $(x+y)^{\alpha}\leqslant2^{\alpha-1}(x^{\alpha
}+y^{\alpha})$ for all $x,y\geqslant0$. Combining this fact with
$b+q_{n}(u)\geqslant b\geqslant1$ and (\ref{qnu0}) thus implies that, for
$K>1$ fixed and all $n$ large enough,%
\[
\frac{1}{K4^{b-1}}\int_{\mathcal{K}_{n}^{>\tau}}\left\vert \beta
_{n}(u)\right\vert ^{b+q_{n}(u)}du\leqslant R_{n}^{X}+R_{n}^{Y}+S_{n}%
^{X}+S_{n}^{Y}%
\]
where%
\[
R_{n}^{X}=\int_{\mathcal{K}_{n}^{>\tau}}\left\vert \frac{B_{n}^{X}(u)}%
{h(u)}\right\vert ^{b+q_{n}(u)}du,\quad S_{n}^{X}=\int_{\mathcal{K}_{n}%
^{>\tau}}\left\vert \frac{Z_{n}^{X}(u)}{h(u)}\right\vert ^{b+q_{n}(u)}du.
\]
It remains to prove that for an appropriate choice of $\overline{u}$ and
$\beta$ we have%
\begin{align*}
\underset{n\rightarrow+\infty}{\lim\sup}\,\mathbb{P}\left(  \mathcal{D}%
_{n}\cap\left\{  R_{n}^{X}\geqslant\frac{\lambda\tau^{s_{n}}}{8c_{0}}\right\}
\right)   &  <\frac{\varepsilon}{8},\\
\underset{n\rightarrow+\infty}{\lim\sup}\,\mathbb{P}\left(  \mathcal{D}%
_{n}\cap\left\{  S_{n}^{X}\geqslant\frac{\lambda\tau^{s_{n}}}{8c_{0}}\right\}
\right)   &  <\frac{\varepsilon}{8},
\end{align*}
which ensures by (\ref{IKnbqn}) that $\mathbb{P}\left(  \mathcal{D}_{n}%
\cap\left\{  v_{n}I_{\mathcal{K}_{n}^{>\tau}}\geqslant\lambda/2\right\}
\right)  \leqslant\varepsilon/2$. For short, it is assumed below that
$1/9<\tau^{s_{n}}/8$.$\smallskip$

\noindent\textbf{(iv)} The following integral $T_{n}$ is crucial with respect to
the integrability of the processes $B_{n}^{X}$ and $Z_{n}^{X}$. Let $b^{\prime
}>b$ be so close to $b$ that $0<b^{\prime}(5-b^{\prime})/6<\gamma<1$. Consider the random function $q_{n}(u)$ from (\ref{sn}). For all
$n$ large enough we have $b\leqslant b+q_{n}(u)<b^{\prime}$ hence (\ref{Lb3})
and (\ref{integamma}) entail%
\begin{align*}
T_{n}  &  =\int_{\mathcal{K}_{n}}\left\vert \frac{\sqrt{1-u}}{h(u)}\right\vert
^{b+q_{n}(u)}du\leqslant\int_{\mathcal{K}_{n}}\left\vert \frac{F^{-1}%
(u)}{\sqrt{1-u}}\right\vert ^{b+q_{n}(u)}du\leqslant\int_{\mathcal{K}_{n}%
}\left\vert \frac{F^{-1}(u)}{\sqrt{1-u}}\right\vert ^{b^{\prime}}du\\
&  \leqslant\left[  -\frac{(1-u)^{1-\gamma}}{1-\gamma}\right]  _{\overline{u}%
}^{1-j_{n}/n}\leqslant\frac{(1-\overline{u})^{1-\gamma}}{1-\gamma}.
\end{align*}
\textbf{(v)} On the one hand we have, by Fubini-Tonelli and recalling that
$B_{n}^{X}$ is a standard Brownian bridge and the sequence $s_n$ is defined at (\ref{qnu0}),%
\begin{align*}
\mathbb{E}\left(  R_{n}^{X}\right)   &  \leqslant T_{n}\sup_{u\in
\mathcal{K}_{n}}\mathbb{E}\left(  \left\vert \frac{B_{n}^{X}(u)}{\sqrt
{u(1-u)}}\right\vert ^{b+q_{n}(u)}\right) \\
&  \leqslant T_{n}\sup_{0\leqslant s\leqslant s_{n}}\mathbb{E}\left(
\left\vert \mathcal{N}\left(  0,1\right)  \right\vert ^{b+s}\right)
=T_{n}\mathbb{E}\left(  \left\vert \mathcal{N}\left(  0,1\right)  \right\vert
^{b+s_{n}}\right)  .
\end{align*}
Assuming $n$ so large that $s_{n}<2-b$ we get $\mathbb{E}\left(  R_{n}%
^{X}\right)  /T_{n}<\mathbb{E}(\left\vert \mathcal{N}\left(  0,1\right)
\right\vert ^{2})=1$ then choosing $\overline{u}_0$ such that $(1-\overline
{u}_0)^{1-\gamma}<8(1-\gamma)\lambda/9c_{0}\varepsilon$ yields, for all $\overline{u} \in [\overline{u}_0,1)$,%
\begin{equation}
\mathbb{P}\left(  R_{n}^{X}\geqslant\frac{\lambda}{9c_{0}}\right)
<\frac{9c_{0}}{\lambda}T_{n}<\frac{\varepsilon}{8}. \label{RXn}%
\end{equation}
On the other hand we have $\mathcal{K}_{n}^{>\tau}\subset\mathcal{K}%
_{n}\subset\Delta_{n}^{\prime}$ and%
\[
S_{n}^{X}\leqslant\sup_{u\in\mathcal{K}_{n}^{>\tau}}\left\vert \frac{Z_{n}%
^{X}(u)}{\sqrt{1-u}}\right\vert ^{b+q_{n}(u)}T_{n}.
\]
By Lemma \ref{Lem_E_approxim} it almost surely holds, for $b^{\prime}>b$ and
all $n$ large,%
\[
\sup_{u\in\mathcal{K}_{n}^{>\tau}}\left\vert \frac{Z_{n}^{X}(u)}{\sqrt{1-u}%
}\right\vert \leqslant\sup_{u\in\mathcal{K}_{n}}\frac{1}{n^{\upsilon}%
\sqrt{1-u}}=\frac{1}{n^{\upsilon}}\sqrt{\frac{n}{j_{n}}}\leqslant
n^{(1-\beta)/2-\upsilon}%
\]
which vanishes provided $1-2\upsilon<\beta<1$. Therefore, for this choice of
$\beta$,%
\begin{equation}
\lim_{n\rightarrow+\infty}S_{n}^{X}=0\ \ a.s.,\quad\lim_{n\rightarrow+\infty
}\mathbb{P}\left(  S_{n}^{X}\geqslant\frac{\lambda}{9c_{0}}\right)  =0.
\label{SXn}%
\end{equation}
\textbf{(vi)} Putting together the conclusions of (i)-(v), and especially (\ref{IKnbqn}%
), (\ref{RXn}) and (\ref{SXn}), implies%
\[
\mathbb{P}\left(  \mathcal{C}_{n}^{\lambda}\right)  \leqslant1-\mathbb{P}%
\left(  \mathcal{D}_{n}\right)  +\mathbb{P}\left(  \mathcal{D}_{n}\cap\left\{
v_{n}I_{\mathcal{K}_{n}^{>\tau}}\geqslant\frac{\lambda}{2}\right\}  \right)
<\frac{\varepsilon}{2}+4\frac{\varepsilon}{8}=\varepsilon.
\]
Finally notice that the same $\beta$ works whatever the choice of $\lambda
,\varepsilon$.$\quad\square\smallskip$

\noindent\textbf{Step 4.}\ Now $\mathcal{L}=\left[  \underline{u},\overline
{u}\right]  $ is fixed. By Lemmas \ref{Lem_tube} and \ref{Lem_E_approxim}
there almost surely exists $n_{0}(\omega)$ such that, for all $n\geqslant
n_{0}(\omega)$, $\varepsilon_{n}(u)$ from (\ref{epsilnu}), $B_n(u)=B_n^X(u)-B_n^Y(u)$ and $Z_n(u)=Z_n^X(u)-Z_n^Y(u)$,
\[
\left\vert \frac{\beta_{n}(u)}{\sqrt{n}}\right\vert \leqslant\varepsilon
_{n}(u)\leqslant x_{0},\quad\beta_{n}(u)=\frac{B_{n}(u)+Z_{n}(u)}{h(u)},\quad
u\in\mathcal{L}.
\]
As a consequence, the cost $\rho_{c}$ is evaluated at $0$ all along this step. Let $\alpha>0$ and consider $I_{\mathcal{L}}=I_{\mathcal{L}_{1,n}%
}+I_{\mathcal{L}_{2,n}}+I_{\mathcal{L}_{3,n}}$ where, for $n\geqslant
n_{0}(\omega)$,%
\begin{equation}
I_{\mathcal{L}_{k,n}}=\int_{\mathcal{L}_{k,n}}\rho_{c}\left(  \frac
{B_{n}(u)+Z_{n}(u)}{\sqrt{n}h(u)}\right)  du,\quad k=1,2,3, \label{ILk}%
\end{equation}
and $\mathcal{L}=\mathcal{L}_{1,n}\cup\mathcal{L}_{2,n}\cup\mathcal{L}_{3,n}$
with $\mathcal{L}_{1,n}=\mathcal{L}\cap\left\{  \left\vert B_{n}(u)\right\vert
\leqslant\alpha\right\}  $, $\mathcal{L}_{2,n}=\mathcal{L}\cap\left\{
\left\vert B_{n}(u)\right\vert \geqslant1/\alpha\right\}  $ and $\mathcal{L}%
_{3,n}=\mathcal{L}\cap\left\{  \alpha<\left\vert B_{n}(u)\right\vert
<1/\alpha\right\}  $. Also define 
$$0<\underline{h}=\min_{u\in\mathcal{L}}h(u)\leqslant \overline{h}=\max_{u\in\mathcal{L}}h(u)<+\infty.$$

\noindent\textbf{Step 4.1} Choose $\alpha\in\left(  0,1\right)  $ arbitrarily
small. In view of the almost sure rate $1/n^{\upsilon}$ from Lemma
\ref{Lem_E_approxim} and (\ref{L}) we have, given $\underline{u}$,
$\overline{u}$ then $\underline{h}$,
\begin{align}
\lim_{n\rightarrow+\infty}v_{n}I_{\mathcal{L}_{1,n}}  &  \leqslant
\lim_{n\rightarrow+\infty}\frac{1}{\rho(1/\sqrt{n})}\int_{\mathcal{L}_{1,n}%
}\rho\left(  \frac{\alpha+1/n^{\upsilon}}{\sqrt{n}\underline{h}}\right)
du\nonumber\\
&  \leqslant\lim_{n\rightarrow+\infty}\frac{\rho(2\alpha /\sqrt{n}\underline{h})}{\rho(1/\sqrt{n})} =\frac{(2\alpha)^{b}}{\underline{h}^{b}}\quad a.s. \label{limvnL1}%
\end{align}
The last equality holds by definition of $\rho\in RV(0,b)$.$\smallskip$

\noindent\textbf{Step 4.2} Write $\mathcal{L}_{2,n}^{+}=\mathcal{L}\cap\left\{
B_{n}(u)\geqslant1/\alpha\right\}  $ and $\mathcal{L}_{2,n}^{-}=\mathcal{L}%
\cap\left\{  B_{n}(u)\leqslant-1/\alpha\right\}  $. By Lemma
\ref{Lem_E_approxim} we have, for $n$ large enough,%
\[
v_{n}I_{\mathcal{L}_{2,n}^{+}}=\frac{1}{\rho(1/\sqrt{n})}\int_{\mathcal{L}%
_{2,n}^{+}}\rho_{+}\left(  \frac{\beta_{n}(u)}{\sqrt{n}}\right)
du\leqslant\frac{1}{\rho(1/\sqrt{n})}\int_{\mathcal{L}_{2,n}^{+}}\rho
_{+}\left(  \frac{2B_{n}(u)}{h(u)\sqrt{n}}\right)  du
\]
then similar arguments as for (ii) in the proof of Lemma \ref{Lem_E_step3}
yield%
\[
v_{n}I_{\mathcal{L}_{2,n}^{+}}\leqslant c_{0}\frac{\rho_{+}(1/\sqrt{n})}%
{\rho(1/\sqrt{n})}\int_{\mathcal{L}_{2,n}^{+}}\left(  \frac{2B_{n}(u)}%
{h(u)}\right)  ^{b_{+}+s_{n}}du
\]
where $s_{n}\rightarrow0$ is defined at (\ref{sn}) with $\tau=2/\alpha$. By
replacing $\min(u,1-u)$ with $u(1-u)\leqslant \min(u,1-u)$ in $(CFG3)$ it
follows that%
\[
v_{n}I_{\mathcal{L}_{2,n}^{+}}\leqslant K\int_{\mathcal{L}}1_{\left\{
B_{n}(u)\geqslant1/\alpha\right\}  }\left\vert \frac{F^{-1}(u)}{\sqrt{u(1-u)}%
}\right\vert ^{b_{+}+s_{n}}\left(  \frac{B_{n}(u)}{\sqrt{u(1-u)}}\right)
^{b_{+}+s_{n}}du
\]
where $K>0$. As a consequence of $(CFG_{E})$ we obtain exactly as for
(\ref{Lb3}) and (\ref{integamma}) that if $b^{\prime}\in(b,2)$ is chosen
sufficiently close to $b$ then
\begin{equation}
\int_{\left(  0,1\right)  }\left\vert \frac{F^{-1}(u)}{\sqrt{u(1-u)}%
}\right\vert ^{b^{\prime}}du=K^{\prime}<+\infty. \label{integb'}%
\end{equation}
Since $2u-1 \leqslant H(u,u) \leqslant u$ for $u \in (0,1)$ we have 
$$-(1-u)^2 \leqslant Cov(B_{n}^X(u),B_n^Y(u))=H(u,u)-u^2 \leqslant u(1-u)$$
hence
$$ 0 \leqslant Var(B_{n}(u)) \leqslant 2u(1-u)+2(1-u)^2=2(1-u) $$
and the $r.v.$ $ B_{n}(u)/\sqrt{u(1-u)} $ is centered Gaussian with variance bounded above by $ 2/\underline{u} $. Let denote $\mathcal{N}(0,1)$ the standard normal distribution. By H\"{o}lder inequality we have, for $u\in\mathcal{L}$ and $n$ large,
\[
\mathbb{E}\left(  1_{\left\{  B_{n}(u)\geqslant1/\alpha\right\}  }\left\vert
\frac{B_{n}(u)}{\sqrt{u(1-u)}}\right\vert ^{b_{+}+s_{n}}\right)
\leqslant K^{\prime \prime}\mathbb{P}\left(  \sup_{\underline{u}\leqslant
u\leqslant\overline{u}}\left\vert B_{n}(u)\right\vert \geqslant\frac{1}%
{\alpha}\right)  ^{1/2}%
\]
where $K^{\prime \prime}=(3/\underline{u})\sup_{b_{+}\leqslant s\leqslant b} \left(\mathbb{E}
\left\vert \mathcal{N}(0,1)\right\vert ^{2s} \right) ^{1/2}<+\infty$
only depends on $b$. We conclude that it asymptotically holds%
\begin{equation}
\mathbb{E}\left(  v_{n}I_{\mathcal{L}_{2,n}^{+}}\right)  \leqslant KK^{\prime
}K^{\prime\prime}\mathbb{P}\left(  \sup_{\underline{u}\leqslant u\leqslant
\overline{u}}\left\vert B_{n}(u)\right\vert >\frac{1}{\alpha}\right)
^{1/2}\leqslant C\exp\left(  -\frac{1}{\alpha^{2}}\right)  \label{EvnL2}%
\end{equation}
where $C$ depends on $M,b,F$ and $\alpha$ was left arbitrary from the beginning. Clearly $\mathbb{E(}%
v_{n}I_{\mathcal{L}_{2,n}^{-}})$ also obeys (\ref{EvnL2}) by the same arguments. Notice that for the left hand tail $\underline{u}$ and $1-\overline{u}$ play a symetric role in the previous control of the variance of $B_n(u)$ by $u(1-u)$.$\smallskip$

\noindent\textbf{Step 4.3} Let introduce $\mathcal{L}_{3,n}^{-}=\mathcal{L}\cap\left\{  -1/\alpha<B_{n}(u)<-\alpha\right\}$ and $\mathcal{L}_{3,n}^{+}=\mathcal{L}\cap\left\{\alpha<B_{n}(u)<1/\alpha\right\}$. By Lemma
\ref{Lem_E_approxim} we almost surely ultimately have
\[
sign(B_{n}(u)+Z_{n}(u))1_{\mathcal{L}_{3,n}}(u)=sign(B_{n}(u))1_{\mathcal{L}%
_{3,n}}(u)
\]
where $sign(x)=1_{x>0}-1_{x<0}$. Therefore, $(C2)$ implies, for all $n$ large
enough,%
\begin{align*}
&  1_{\mathcal{L}_{3,n}}(u)\rho_{c}\left(  \frac{B_{n}(u)+Z_{n}(u)}{\sqrt
{n}h(u)}\right) \\
&  =1_{\mathcal{L}_{3,n}^{+}}(u)\rho_{+}\left(  \frac{B_{n}(u)+Z_{n}(u)}%
{\sqrt{n}h(u)}\right)  +1_{\mathcal{L}_{3,n}^{-}}(u)\rho_{-}\left(
\frac{\left\vert B_{n}(u)+Z_{n}(u)\right\vert }{\sqrt{n}h(u)}\right)  .
\end{align*}
Now assume that $\alpha<2/\overline{h}$ and $\mathcal{L}_{3,n}\neq\emptyset$, so
that%
\begin{align*}
v_{n}I_{\mathcal{L}_{3,n}}  &  =\frac{1}{\rho(1/\sqrt{n})}\left(
\int_{\mathcal{L}_{3,n}^{+}}\rho_{+}\left(  \frac{\left\vert B_{n}%
(u)\right\vert }{\sqrt{n}h(u)}\right)  du+R_{n}^{+}\right) \\
&  +\frac{1}{\rho(1/\sqrt{n})}\left(  \int_{\mathcal{L}_{3,n}^{-}}\rho
_{-}\left(  \frac{\left\vert B_{n}(u)\right\vert }{\sqrt{n}h(u)}\right)
du+R_{n}^{-}\right)
\end{align*}
where we have, by convexity and differentiability of $\rho_{\pm}$ on
$(0,+\infty)$,
\begin{align*}
R_{n}^{\pm}  &  =\int_{\mathcal{L}_{3,n}^{\pm}}\left(  \rho_{\pm}\left(
\frac{\left\vert B_{n}(u)+Z_{n}(u)\right\vert }{\sqrt{n}h(u)}\right)
-\rho_{\pm}\left(  \frac{\left\vert B_{n}(u)\right\vert }{\sqrt{n}%
h(u)}\right)  \right)  du\\
&  \leqslant\sup_{u\in\mathcal{L}_{3,n}^{\pm}}\rho_{\pm}^{\prime}\left(
\frac{\left\vert B_{n}(u)\right\vert +\left\vert Z_{n}(u)\right\vert }%
{\sqrt{n}h(u)}\right)  \frac{\left\vert Z_{n}(u)\right\vert }{\sqrt{n}h(u)}%
\end{align*}
The regular variation $(C2)$ further implies $x\rho_{\pm}^{\prime}%
(x)/\rho_{\pm}(x)\rightarrow1$ as $x\rightarrow0$. As a consequence, with
probability one, for all $n$ large it holds
\begin{align*}
\frac{R_{n}^{\pm}}{\rho(1/\sqrt{n})}  &  \leqslant\frac{1}{\rho(1/\sqrt{n}%
)}\rho_{\pm}\left(  \frac{\left\vert B_{n}(u)\right\vert +\left\vert
Z_{n}(u)\right\vert }{\sqrt{n}h(u)}\right)  \sup_{u\in\mathcal{L}_{3,n}}%
\frac{\left\vert Z_{n}(u)\right\vert }{\left\vert B_{n}(u)\right\vert
+\left\vert Z_{n}(u)\right\vert }\\
&  \leqslant\frac{\rho_{\pm}\left(  2/\sqrt{n}\underline{h}\alpha\right)
}{\rho(1/\sqrt{n})}\frac{2}{\alpha n^{\upsilon}}\leqslant\frac{\rho_{\pm
}\left(  2/\sqrt{n}\underline{h}\alpha\right)  }{\rho_{\pm}(1/\sqrt{n})}%
\frac{2}{\alpha n^{\upsilon}}\leqslant\left(  \frac{2}{\underline{h}\alpha
}\right)  ^{b_{\pm}}\frac{3}{\alpha n^{\upsilon}}%
\end{align*}
which vanishes as $n\rightarrow+\infty$. Here we have used that $\rho_{\pm
}\left(  \theta x\right)  /\rho_{\pm}\left(  x\right)  \rightarrow
\theta^{b_{\pm}}$ as $x\rightarrow0$ for any fixed $\theta>0$, and Lemma
\ref{Lem_E_approxim}. Finally we see that%
\[
\frac{1}{\rho_{\pm}(1/\sqrt{n})}\int_{\mathcal{L}_{3,n}^{\pm}}\rho_{\pm
}\left(  \frac{\left\vert B_{n}(u)\right\vert }{\sqrt{n}h(u)}\right)
du=\int_{\mathcal{L}_{3,n}^{\pm}}\left(  \frac{\left\vert B_{n}(u)\right\vert
}{h(u)}\right)  ^{b}du+R_{3,n}^{\pm}%
\]
with%
\[
R_{3,n}^{\pm}=\int_{\mathcal{L}_{3,n}^{\pm}}L_{n}^{\pm}(u)\left(
\frac{\left\vert B_{n}(u)\right\vert }{h(u)}\right)  ^{b}du,\quad L_{n}^{\pm
}(u)=\frac{L_{\pm}(\left\vert B_{n}(u)\right\vert /\sqrt{n}h(u))}{L_{\pm
}(1/\sqrt{n})}-1.
\]
Clearly, it follows
\[
\left\vert R_{3,n}^{\pm}\right\vert \leqslant\left(  \frac{1}{\underline
{h}\alpha}\right)  ^{b}\sup_{u\in\mathcal{L}_{3,n}}\left\vert L_{n}^{\pm
}(u)\right\vert \leqslant\left(  \frac{1}{\underline{h}\alpha}\right)
^{b}\left(  \sup_{\alpha/\overline{h}\sqrt{n}\leqslant x\leqslant
1/\alpha\underline{h}\sqrt{n}}\frac{L_{\pm}(x)}{L_{\pm}(1/\sqrt{n})}-1\right)
\]
thus, by (\ref{varlente}) and (\ref{eta0}) we get $\left\vert R_{3,n}^{\pm}\right\vert \rightarrow0$
as $n\rightarrow+\infty$. We conclude that%
\begin{equation}
I_{\mathcal{L}_{3,n}}^{\ast}=\rho_{+}(1/\sqrt{n})\int_{\mathcal{L}_{3,n}^{+}%
}\left(  \frac{\left\vert B_{n}(u)\right\vert }{h(u)}\right)  ^{b_{+}}%
du+\rho_{-}(1/\sqrt{n})\int_{\mathcal{L}_{3,n}^{-}}\left(  \frac{\left\vert
B_{n}(u)\right\vert }{h(u)}\right)  ^{b_{-}}du \label{L3nstar}%
\end{equation}
almost surely satisfies $\lim_{n\rightarrow+\infty}v_{n}\vert I_{\mathcal{L}_{3,n}%
}-I_{\mathcal{L}_{3,n}}^{\ast}\vert =0.\smallskip$

\noindent\textbf{Step 5.} Consider $W_{c}(\mathbb{F}_{n},\mathbb{G}%
_{n})=\int_{0}^{1}\rho_{c}(\mathbb{F}_{n}^{-1}(u),\mathbb{G}_{n}^{-1}(u))du$.
As we assumed that%
\[
\lim_{n\rightarrow+\infty} \frac{\rho_{+}(1/\sqrt{n})}{\rho(1/\sqrt{n})}=\pi_{+},\quad
\lim_{n\rightarrow+\infty}\frac{\rho_{-}(1/\sqrt{n})}{\rho(1/\sqrt{n})}=\pi_{-}%
\]
by $(C4)$ and $E=\mathbb{R}$ we have established that $v_{n}W_{c}^{E}(\mathbb{F}%
_{n},\mathbb{G}_{n})\rightarrow_{weak}W$ with%
\[
W=\pi_{+}\int_{0}^{1}1_{\left\{  B(u)>0\right\}  }\left(  \frac{\left\vert
B(u)\right\vert }{h(u)}\right)  ^{b_{+}}du+\pi_{-}\int_{0}^{1}1_{\left\{
B(u)<0\right\}  }\left(  \frac{\left\vert B(u)\right\vert }{h(u)}\right)
^{b_{-}}du
\]
and $B$ is a standard Brownian bridge. To see this write $W_{c}^{E}%
(\mathbb{F}_{n},\mathbb{G}_{n})=I_{\mathcal{I}_{n}}+I_{\mathcal{J}_{n}%
}+I_{\mathcal{K}_{n}}+I_{\mathcal{L}_{1,n}}+I_{\mathcal{L}_{2,n}%
}+I_{\mathcal{L}_{3,n}}$ where each of the first three integrals is indeed the
sum of its left hand tail and right hand tail version, likewise for $I_{\mathcal{L}_{3,n}%
}^{\ast}$ defined at (\ref{L3nstar}). We have shown that $v_{n}(I_{\mathcal{I}%
_{n}}+I_{\mathcal{J}_{n}})\rightarrow0$ in probability.\ Let $\Psi$ be a real
valued $k$-Lipschitz function on $\mathbb{R}$, bounded by $m$. Given
arbitrarily small constants $\lambda>0$, $\varepsilon>0$ and $\alpha>0$ then
an appropriate choice of $0<\underline{u},\overline{u}<1$ and thus
$\underline{h}$ it holds, for all $n$ large enough, by Lemma \ref{Lem_E_step3} and Step 4,
\begin{align*}
&  \mathbb{E}\left(  \Psi\left(  v_{n}(I_{\mathcal{K}_{n}}+I_{\mathcal{L}%
_{1,n}}+I_{\mathcal{L}_{2,n}}+I_{\mathcal{L}_{3,n}})\right)  -\Psi\left(
v_{n}I_{\mathcal{L}_{3,n}}^{\ast}\right)  \right)  \\
&  \leqslant4m\mathbb{P}\left(  v_{n}I_{\mathcal{K}_{n}}>\lambda\right)
+4m\mathbb{P}\left(  v_{n}\left\vert I_{\mathcal{L}_{3,n}}-I_{\mathcal{L}%
_{3,n}}^{\ast}\right\vert >\lambda\right)  \\
&  +4m\mathbb{P}\left(  v_{n}I_{\mathcal{L}_{1,n}}>\frac{(5\alpha)^{b}%
}{\underline{h}^{b}}\right)  +k\mathbb{E}\left(  4\lambda+\frac{(5\alpha)^{b}%
}{\underline{h}^{b}}+v_{n}I_{\mathcal{L}_{2,n}}\right)  \\
&  \leqslant12m\varepsilon+4k\lambda+\frac{k(5\alpha)^{b}}{\underline{h}^{b}%
}+2kC\exp\left(  -\frac{1}{\alpha^{2}}\right)
\end{align*}
which is as small as desired. Finally it is easilly seen that $v_{n}%
I_{\mathcal{L}_{3,n}}^{\ast}\rightarrow_{weak}W$ as $(\underline{u}%
,\overline{u})\rightarrow(0,1)$ and $\alpha\rightarrow0$ so that
$\mathbb{E(}\Psi\left(  W\right)  )$ can replace $\mathbb{E(}\Psi
(v_{n}I_{\mathcal{L}_{3,n}}^{\ast}))$ above with an asymptotically arbitrarily
small error.$\quad\square$

\subsection{The case $F<G$}\label{proofFneG}
We establish Theorem \ref{D=R}.
\smallskip 

\noindent\textbf{Step 0.} In this section $D=(0,1)$. Without loss of generality, assume that $F^{-1}>G^{-1}$ everywhere. We again focus on arguments for the right hand tail, thus we write $\psi _X=\psi _X^+$ and $\psi _Y=\psi _Y^+$ on $(y_0,+\infty)$. Therefore $\psi _X^{-1} > \psi _Y^{-1}$ and $\psi _X^{-1} > 0$ on $(u_0,1)$ where $u_{0}=F^{-1}(y_0)$. We need this stochastic ordering only to simplify the control of extremes without imposing $(CFG_{E})$. Let assume $(FG)$, $(C)$ with $b\in\left[  1,2\right)  $ and $(CFG_{D})$. For $y$ large it holds $\rho_{\pm}\left(  y\right)  =\exp(l_{\pm}(y))$ with $l_{\pm}\in{RV}_{2}^{+}(\gamma_{\pm},+\infty)$. By (\ref{CFGD_derivee}), for $y_{0}>0$ and $\theta_{+},\theta_{-}>1$ playing exactly the role of $\theta$ in (CFG) of \cite{BFK17} we have
\begin{equation}
(\psi _X\circ l_{+}^{-1})^{\prime}(y)\geqslant2+\frac{2\theta_{+}}{y},\quad
(\psi _Y\circ l_{-}^{-1})^{\prime}(y)\geqslant2+\frac{2\theta_{-}}{y},\quad
y>y_{0}.\label{CFGDproof}%
\end{equation}
In particular, this implies
\begin{equation}
l_+\circ\psi_X^{-1}(y)\leqslant\frac{y}{2}-\theta_+\log y+K,\quad y>y_{0}.\label{CFGDproof1}%
\end{equation}
By (\ref{CFGD_integree}), whenever $F^{-1}(u)-G^{-1}(u)>0$ is not asymptotically away from $0$ as $u \rightarrow 1$ we further ask that, for some $\theta_2>0$,
\begin{equation}
l_+\circ\psi_X^{-1}(y)\leqslant\frac{y}{2}-2\log\psi_X^{-1}(y)-\theta_2\log y,\quad y>y_{0}.\label{CFGDproof2}
\end{equation}
Notice that if $F$ is logconvex then $\log\psi_X^{-1}(y)>\log y$ and (\ref{CFGDproof2}) already implies (\ref{CFGDproof1}) with $\theta_+>2$ whereas if $F$ is logconcave then $\log\psi_X^{-1}(y)<\log y$ and (\ref{CFGDproof1}) implies (\ref{CFGDproof2}) with $\theta_2>1$. Since $(CFG_D)$ implies $(CFG)$ of \cite{BFK17} through (\ref{CFGD_integree}) hence (\ref{CFGDproof}), we are allowed to use most results of the latter paper. In particular Theorem \ref{D=R} is true when $F^{-1}(u)-G^{-1}(u)>\delta$ for some $\delta>0$ and $b>1$ to ensure $(C3)$ in \cite{BFK17}. We thus focus on the case $F^{-1}(u)-G^{-1}(u) \rightarrow 0$ as $u \rightarrow 1$ which requires (\ref{CFGDproof2}) whatever $b$, and we isolate out the case $b=1$ only when necessary to extend the main result of \cite{BFK17}, at Step 4. We often use $F^{-1}(u)=\psi_{X}^{-1}(\log(1/(1-u)))$. A consequence is that (\ref{CFGDproof2}) also reads
\[
\rho\circ  F^{-1}(u)=\rho_+\circ F^{-1}(u)\leqslant\frac{1}{F^{-1}(u)^{2}\sqrt{1-u}%
\left\vert \log(1-u)\right\vert ^{\theta_2}},\quad u>u_{0}.
\]
Let us study $W_{c}(\mathbb{F}_{n},\mathbb{G}_{n})-W_{c}(F,G)=I_{\mathcal{I}_{n}}+I_{\mathcal{J}_{n}}+I_{\mathcal{K}_{n}}+I_{\mathcal{L}}$ with the notation
\begin{align}
I_{A}  &  =\int_{A}\left(  \rho_{c}\left(  \tau(u)+\tau_{n}(u)\right)-\rho_{c}\left(  \tau(u)\right)  \right)du, \quad A\subset\left(0,1\right), \label{IA} \\
\tau(u)  &  =F^{-1}(u)-G^{-1}(u),\quad\tau_{n}(u)=\frac{\beta_{n}(u)}{\sqrt{n}}=\frac{\beta_{n}^{X}(u)-\beta_{n}^{Y}(u)}{\sqrt{n}},\nonumber\
\end{align}
and $\mathcal{I}_{n}=\left(  1-i_{n}/n,1\right]  $, $\mathcal{J}_{n}=\left(
1-j_{n}/n,1-i_{n}/n\right]  $, $\mathcal{K}_{n}=\left(  \overline{u}%
,1-j_{n}/n\right]  $, $\mathcal{L}=\left[  \underline{u},\overline{u}\right]
$ with $0<\underline{u}<1/2<\overline{u}<1$.$\smallskip$

\noindent\textbf{Step 1.} Consider a non negative increasing sequence $K_{n}\rightarrow+\infty$ to be chosen later in such a way that $K_{n}/\log\log n\rightarrow0$. Define%
\begin{equation}
i_{n}=\frac{\sqrt{n}}{K_{n}\exp\left(  l\circ\psi_{X}^{-1}(\log n+K_{n}%
)\right)  }. \label{inbis}%
\end{equation}
We have $l\circ\psi_{X}^{-1}(y)= l_+\circ\psi_{X}^{-1}(y)\rightarrow+\infty$ as $y\rightarrow+\infty$ thus $i_{n}=o\left(  \sqrt{n}/K_{n}\right)  $. When (\ref{CFGDproof2}) is enforced then for any $\theta^{\prime}\in\left(  0,\theta_2\right)$ and all $n$ large enough,
\begin{align}
i_{n}  &  \geqslant\frac{K}{K_{n}}\left(  F^{-1}\left(  1-\frac{1}{ne^{K_{n}}%
}\right)  \right)  ^{2}\exp\left(  -\frac{K_{n}}{2}+\theta_+\log(\log
n+K_{n})\right) \nonumber\\
&  >\left(  F^{-1}\left(  1-\frac{1}{n}\right)  \right)  ^{2}(\log
n)^{\theta^{\prime}}. \label{inbibis}%
\end{align}
Otherwise, when only (\ref{CFGDproof1}) holds then for $\theta^{\prime}\in\left(  1,\theta_+\right)$, 
\begin{equation}
i_{n} \geqslant\frac{K}{K_{n}}\exp\left(  -\frac{K_{n}}{2}+\theta_+\log(\log(n+K_{n})) \right)>(\log n)^{\theta^{\prime}}. \label{inbibis2}%
\end{equation}
Hence in both case we have ${i_{n}/\log\log n\rightarrow+\infty}$ and $i_{n}/\sqrt
{n}\rightarrow0$. Let us define
\begin{align*}
I_{\mathcal{I}_{n}}^{1}  &  =\int_{\mathcal{I}_{n}}\rho_{+}\left(
\tau(u)\right)  du,\\
I_{\mathcal{I}_{n}}^{2}  &  =\int_{\mathcal{I}_{n}}\rho_{c}\left(
\mathbb{F}_{n}^{-1}(u)-\mathbb{G}_{n}^{-1}(u)\right)  du=\frac{1}{n}%
{\sum\limits_{i=n-[i_{n}]}^{n}}\rho_{c}\left(  X_{(i)}-Y_{(i)}\right)  .
\end{align*}

\begin{lemma}
\label{lem:DS}Assume that $(C),$ $(FG)$ and $(CFG_{D})$ hold. Then $\sqrt
{n}I_{\mathcal{I}_{n}}^{1}\rightarrow0$ and $\sqrt{n}I_{\mathcal{I}_{n}}%
^{2}\rightarrow0$ in probability.
\end{lemma}

\noindent\textbf{Proof.} This readily follows from Lemma 22 in \cite{BFK17}. For $\sqrt{n}I_{\mathcal{I}_{n}}^{1}$ the mentioned proof
only needed $\theta>0$ hence $\theta_+,\theta_->0$. For $\sqrt{n}I_{\mathcal{I}_{n}}^{2}$ the initial
expansion%
\[
{\sum\limits_{i=n-[i_{n}]}^{n}}\rho_{c}\left(  X_{(i)}-Y_{(i)}\right)
\leqslant{\sum\limits_{i=n-[i_{n}]}^{n}}\rho_{+}\left(  X_{(i)}\right)
+{\sum\limits_{i=n-[i_{n}]}^{n}}\rho_{-}\left(  Y_{(i)}\right)
\]
almost surely holds for $n$ large enough, when $\min(X_{(n-[i_{n}%
])},Y_{(n-[i_{n}])})>0$.$\quad\square\smallskip$

\noindent\textbf{Step 2}. We now study $I_{\mathcal{J}_{n}}$ with
$j_{n}=n^{\beta}$, $\beta\in\left(  1/2,1\right)  $. Recall that $\Delta
_{n}=\mathcal{J}_{n}\cup\mathcal{K}_{n}\cup\mathcal{L}$ and $\tau
(u)=F^{-1}(u)-G^{-1}(u)>0$ for all $u\in\Delta_{n}$. $\smallskip$

\noindent\textbf{(i)} Define $\varepsilon_{n}=\sup_{u\in\Delta_{n}%
}\varepsilon_{n}(u)$ where $\varepsilon_{n}(u)=\varepsilon_{n}^{X}%
(u)+\varepsilon_{n}^{Y}(u)$ and%
\[
\varepsilon_{n}^{X}(u)=\frac{\sqrt{\log\log n}}{\sqrt{n}}\frac{\sqrt{1-u}%
}{h_{X}(u)},\quad\varepsilon_{n}^{Y}(u)=\frac{\sqrt{\log\log n}}{\sqrt{n}%
}\frac{\sqrt{1-u}}{h_{Y}(u)}.
\]
The current $\varepsilon_{n}$ is bounded by the one of (\ref{epsiln}). By combining (\ref{hongr}) and (\ref{epsiln}) with (\ref{inbibis}) as in Lemma
\ref{Lem_tube} we get, for some $\zeta>0$,%
\[
\lim_{n\rightarrow+\infty}(\log n)^{\zeta}\sup_{u\in\Delta_{n}}\tau
_{n}(u)\leqslant9\lim_{n\rightarrow+\infty}(\log n)^{\zeta}\varepsilon
_{n}=0\quad a.s.
\]
Let $m_{n}\rightarrow+\infty$ be a non negative sequence so slow that
$m_{n}\varepsilon_{n}\rightarrow0$. Consider $\mathcal{J}_{n}=\mathcal{J}%
_{n}^{<}\cup\mathcal{J}_{n}^{>}$ where%
\begin{align*}
\mathcal{J}_{n}^{<} &  =\left\{  u\in\mathcal{J}_{n}:0<\tau(u)\leqslant
m_{n}\varepsilon_{n}(u)\right\}  ,\\
\mathcal{J}_{n}^{>} &  =\left\{  u\in\mathcal{J}_{n}:0<m_{n}\varepsilon
_{n}(u)<\tau(u)\right\}  .
\end{align*}
By (\ref{hongr}) again we almost surely ultimately have
\[
-9\varepsilon_{n}(u)<\tau_{n}(u)=\frac{\beta_{n}^{X}(u)}{\sqrt{n}}-\frac
{\beta_{n}^{Y}(u)}{\sqrt{n}}<9\varepsilon_{n}(u),\quad u\in\mathcal{J}_{n}.
\]
Notice that if $u\in\mathcal{J}_{n}^{>}$ then%
\begin{equation}
0<(m_{n}-9)\varepsilon_{n}(u)<\tau(u)+\tau_{n}(u)<\tau(u)+9\varepsilon
_{n}(u)<\tau(u)\left(  1+\frac{9}{m_{n}}\right)  \label{DFnplus}%
\end{equation}
whereas if $u\in\mathcal{J}_{n}^{<}$ then it is possible that $\tau
(u)+\tau_{n}(u)<0$ since%
\begin{equation}
-9\varepsilon_{n}(u)<\tau_{n}(u)<\tau(u)+\tau_{n}(u)<(m_{n}+9)\varepsilon
_{n}(u).\label{DFnmoins}%
\end{equation}
Let us control $\left\vert I_{\mathcal{J}_{n}}\right\vert \leqslant\left\vert
I_{\mathcal{J}_{n}^{<}}\right\vert +\left\vert I_{\mathcal{J}_{n}^{>}%
}\right\vert $, starting with the first term.$\smallskip$

\noindent\textbf{(ii)} Recall that $\sup_{u\in\mathcal{J}_{n}^{<}}%
m_{n}\varepsilon_{n}(u)\rightarrow0$ as $n\rightarrow+\infty$. By
(\ref{DFnmoins}) we have, for $u\in I_{\mathcal{J}_{n}^{<}}$ and $m_{n}>9$,%
\begin{align*}
\left\vert \rho_{c}\left(  \tau(u)+\tau_{n}(u)\right)  -\rho_{c}\left(
\tau(u)\right)  \right\vert  &  \leqslant\rho_{c}\left(  \tau(u)+\tau
_{n}(u)\right)  +\rho_{+}\left(  \tau(u)\right)  \\
&  \leqslant\rho_{-}\left(  9\varepsilon_{n}(u)\right)  +2\rho_{+}%
(2m_{n}\varepsilon_{n}(u)).
\end{align*}
hence $\sqrt{n}\left\vert I_{\mathcal{J}_{n}^{<}}\right\vert \leqslant
R_{1,n}+R_{2,n}$ for all $n$ large enough, with%
\begin{align*}
R_{1,n} &  =K\sqrt{n}\int_{\mathcal{J}_{n}^{<}}\varepsilon_{n}(u)^{b_{-}}%
L_{-}\left(  9\varepsilon_{n}(u)\right)  du,\\
R_{2,n} &  =K\sqrt{n}\int_{\mathcal{J}_{n}^{<}}(m_{n}\varepsilon
_{n}(u))^{b_{+}}L_{+}\left(  2m_{n}\varepsilon_{n}(u)\right)  du.
\end{align*}

\begin{lemma}
\label{Lem_D_step2}Assume $(C)$, $(FG)$ and $(CFG)$. We have $R_{1,n}\rightarrow0$ and $R_{2,n}\rightarrow0$.
\end{lemma}

\noindent\textbf{Proof.} If $F^{-1}(u)-G^{-1}(u)>\delta$ then the set $\mathcal{J}_{n}^{<}$ is ultimately empty. Otherwise (\ref{CFGDproof2}) holds. We have $\sqrt{1-u}\left(1/h_{X}(u)+1/h_{Y}(u)\right)  \leqslant2F^{-1}(u)/\sqrt{1-u}$ for $u\in\mathcal{J}_{n}$ in view of $F^{-1}(u)>G^{-1}(u)$ and $(FG3)$. If $\min(b_{+},b_{-})-1>0$ this extra power cancels the slowly varying functions and we asymptotically have 
\[
R_{1,n}+R_{2,n}\leqslant K\sqrt{n}\int_{\mathcal{J}_{n}^{<}}m_{n}%
\varepsilon_{n}(u)du\leqslant Km_{n}\sqrt{\log\log n}\int_{\mathcal{J}_{n}%
}\frac{F^{-1}(u)}{\sqrt{1-u}}du.
\]
If $b_{+}=1$ then $L_{+}(x)$ is bounded on $\left[  0,x_{0}\right]  $ since
$xL_{+}(x)$ is convex non negative and starts from $0$. Hence $L_{+}\left(
2m_{n}\varepsilon_{n}(u)\right)  $ is bounded on $\mathcal{J}_{n}$, and the
above upper bound remains true. Likewise if $b_{-}=1$ then $L_{-}\left(
9\varepsilon_{n}(u)\right)  $ is bounded on $\mathcal{J}_{n}$. Observe that
(\ref{CFGDproof2}) and $l(y)>\log y$ imply%
\[
\psi^{-1}_X(y)\leqslant\exp(l\circ\psi^{-1}_X(y))\leqslant\frac{1}{\psi
^{-1}_X(y)^{2}}\exp\left(  \frac{y}{2}-\theta\log y\right)
\]
thus $\psi^{-1}_X(y)^{6}\leqslant e^{y}$ and $F^{-1}(u)<1/(1-u)^{1/6}$.
Therefore%
\[
\int_{\mathcal{J}_{n}}\frac{F^{-1}(u)}{\sqrt{1-u}}du\leqslant K\left(
\frac{j_{n}}{n}\right)  ^{1/3}=Kn^{(\beta-1)/3}%
\]
with $\beta<1$ and the conclusion follows since $m_{n}\rightarrow+\infty$ is
arbitrarily slow.$\quad\square$\textit{\smallskip}

We have shown that $\sqrt{n}I_{\mathcal{J}_{n}^{<}}\rightarrow0$
almost surely.$\smallskip$

\noindent\textbf{(iii) }By (\ref{DFnplus}) we ultimately have, for all
$u\in\mathcal{J}_{n}$,%
\[
\left\vert \rho_{c}\left(  \tau(u)+\tau_{n}(u)\right)  -\rho_{c}\left(
\tau(u)\right)  \right\vert =\left\vert \rho_{+}\left(  \tau(u)+\tau
_{n}(u)\right)  -\rho_{+}\left(  \tau(u)\right)  \right\vert .
\]
Consider now $\mathcal{J}_{n}^{>}=\mathcal{J}_{n}^{<\delta}\cup\mathcal{J}%
_{n}^{>\delta}$ with%
\[
\mathcal{J}_{n}^{<\delta}=\left\{  u\in\mathcal{J}_{n}:m_{n}\varepsilon
_{n}(u)<\tau(u)<\delta\right\}  ,\quad\mathcal{J}_{n}^{>\delta}=\left\{
u\in\mathcal{J}_{n}:\tau(u)>\delta\right\}  .
\]
Since $\tau(u)>\delta$ on $\mathcal{J}_{n}^{>\delta}$ and Proposition 31 and
Lemma 25 of \cite{BFK17} are satisfied by $\rho_{+}$ -- thanks to
(\ref{L'}) and (\ref{L1}) -- we readily deduce from Lemma 23 of \cite{BFK17} that
\[\lim_{n\rightarrow+\infty}\sqrt{n}\int_{\mathcal{J}_{n}^{>\delta}}\left\vert
\rho_{+}\left(  \tau(u)+\tau_{n}(u)\right)  -\rho_{+}\left(  \tau(u)\right)
\right\vert du=0\quad a.s.
\]
Concerning $\mathcal{J}_{n}^{<\delta}$ observe that by (\ref{DFnplus}) again
$0<\tau(u)+\tau_{n}(u)<2\delta$ for all $n$ large. Since $\rho_{+}$ is convex
it ensues%
\begin{align*}
\left\vert \rho_{+}\left(  \tau(u)+\tau_{n}(u)\right)  -\rho_{+}\left(
\tau(u)\right)  \right\vert &  \leqslant\max\left(  \rho_{+}^{\prime}\left(  \tau(u)+\tau_{n}(u)\right)
,\rho_{+}^{\prime}\left(  \tau(u)\right)  \right)  \left\vert \tau
_{n}(u)\right\vert \\
&  \leqslant K_{\delta}\left\vert \tau_{n}(u)\right\vert
\end{align*}
with $K_{\delta}=\rho_{+}^{\prime}\left(  2\delta\right)  $. Therefore, with
probability one, for all $n$ large enough
\[
\sup_{u\in\mathcal{J}_{n}^{<\delta}}\left\vert \rho_{+}\left(  \tau
(u)+\tau_{n}(u)\right)  -\rho_{+}\left(  \tau(u)\right)  \right\vert \leqslant
K_{\delta}\sup_{u\in\mathcal{J}_{n}^{<\delta}}\left\vert \tau_{n}(u)\right\vert 
\leqslant K \sup_{u\in\mathcal{J}_{n}^{<\delta}}\left\vert \varepsilon_{n}(u)\right\vert .
\]
As already seen, (\ref{CFGD_integree}) implies $F^{-1}%
(u)<1/(1-u)^{1/6}$ for all $u<1$ large enough. As a consequence, with probability one it ultimately
holds
\begin{align*}
&  \sqrt{n}\int_{\mathcal{J}_{n}^{<\delta}}\left\vert \rho_{+}\left(
\tau(u)+\tau_{n}(u)\right)  -\rho_{+}\left(  \tau(u)\right)  \right\vert
du\leqslant K\sqrt{n}\int_{\mathcal{J}_{n}}\varepsilon_{n}(u)du\\
&  \leqslant K\sqrt{\log\log n}\int_{\mathcal{J}_{n}}\frac{F^{-1}(u)}%
{\sqrt{1-u}}du\leqslant Kn^{(\beta-1)/3}\sqrt{\log\log n}%
\end{align*}
which vanishes as $n\rightarrow+\infty$. We conclude that $\sqrt
{n}I_{\mathcal{J}_{n}^{>}}\rightarrow0$ almost surely.$\smallskip$

\noindent\textbf{Step 3}. The convergence of $I_{\mathcal{K}_{n}}$ is weaker
than in probability.

\begin{lemma}
\label{Lem_D_step3}Assume $(FG)$, $(C)$ and $(CFG_{D})$. There exists
$\beta\in\left(  1/2,1\right)  $ such that for any choice of $\lambda>0$ and
$\varepsilon>0$ one can find $\overline{u}_0\in\left(  1/2,1\right)  $ and
$n_{0}>0$ such that, for all $\overline{u} \in [\overline{u}_0,1)$ and all $n>n_{0}$,
\[
\mathbb{P}\left(  \sqrt{n}I_{\mathcal{K}_{n}}>\lambda\right)  <\varepsilon
\text{.}%
\]

\end{lemma}

\noindent\textbf{Proof}. Fix $\delta>0$ and consider%
\[
\mathcal{K}_{n}^{<\delta}=\left\{  u\in\mathcal{K}_{n}:0<\tau(u)<\delta
\right\}  ,\quad\mathcal{K}_{n}^{>\delta}=\left\{  u\in\mathcal{K}_{n}%
:\tau(u)>\delta\right\}  .
\]
The claimed result holds for $I_{\mathcal{K}_{n}^{>\delta}}$ by applying Lemma 26 from \cite{BFK17} with $\delta=\tau_{0}$ and $\overline
{u}=F(M)$. Let us apply Lemma \ref{Lem_E_approxim} to get, for $K>\sup
_{\left\vert x\right\vert <2\delta}\rho_{c}^{\prime}(x)$,%
\begin{align*}
\sqrt{n}I_{\mathcal{K}_{n}^{<\delta}}  &  =\sqrt{n}\int_{\mathcal{K}%
_{n}^{<\delta}}\left\vert \rho_{c}\left(  \tau(u)+\tau_{n}(u)\right)
-\rho_{c}\left(  \tau(u)\right)  \right\vert du\leqslant K\int_{\mathcal{K}%
_{n}^{<\delta}}\left\vert \beta_{n}(u)\right\vert du\\
&  \leqslant K\int_{\mathcal{K}_{n}}\left(  \frac{\left\vert B_{n}%
^{X}(u)\right\vert }{h_{X}(u)}+\frac{\left\vert B_{n}^{Y}(u)\right\vert
}{h_{Y}(u)}\right)  du+\int_{\mathcal{K}_{n}}\left(  \frac{\left\vert
Z_{n}^{X}(u)\right\vert }{h_{X}(u)}+\frac{\left\vert Z_{n}^{Y}(u)\right\vert
}{h_{Y}(u)}\right)  du.
\end{align*}
The first two terms satisfy%
\begin{align*}
\mathbb{E}\left(  \int_{\mathcal{K}_{n}}\frac{\left\vert B_{n}^{X}%
(u)\right\vert }{h_{X}(u)}du\right)   &  \leqslant\int_{\mathcal{K}_{n}}%
\frac{\sqrt{1-u}}{h_{X}(u)}\frac{\mathbb{E}\left(  \left\vert B_{n}%
^{X}(u)\right\vert \right)  }{\sqrt{u(1-u)}}du\\
&  \leqslant\int_{\overline{u}}^{1}\frac{F^{-1}(1-u)}{\sqrt{1-u}}%
du\leqslant3\left(  1-\overline{u}\right)  ^{1/3}%
\end{align*}
and the last two terms obey, with probability one as $n\rightarrow+\infty$,%
\begin{align*}
\int_{\mathcal{K}_{n}}\frac{\left\vert Z_{n}^{X}(u)\right\vert }{h_{X}(u)}du
&  \leqslant\sup_{u\in\mathcal{K}_{n}}\left\vert Z_{n}^{X}(u)\right\vert
\int_{\mathcal{K}_{n}}\frac{F^{-1}(1-u)}{1-u}du\\
&  \leqslant\frac{1}{n^{\upsilon}}\int_{\overline{u}}^{1-j_{n}/n}\frac
{F^{-1}(1-u)}{1-u}du\\
&  \leqslant\frac{1}{n^{\upsilon}}\int_{\overline{u}}^{1-j_{n}/n}\frac
{1}{(1-u)^{7/6}}du\leqslant\frac{6}{n^{\upsilon}}n^{(1-\beta)/6}%
\end{align*}
which vanishes if $\beta>1-6\upsilon$ is chosen close enough to $1$%
.$\quad\square\medskip$

\noindent\textbf{Step 4}. Here we recall that $(C2)$ with $b_{\pm}>1$ and (\ref{CFGD_derivee}) respectively imply $(C3)$ and $(CFG)$ in \cite{BFK17}. Clearly Steps 4 and 5 of \cite{BFK17} remain true in the current framework and lead to the same conclusion as the main theorem in the latter paper, whence Theorem \ref{D=R}. The new case to conclude with is $b=1$. By Glinvenko-Cantelli and $(FG)$, we almost surely have
\[
0\leqslant\left\vert \tau_{n}(u)\right\vert = \frac{\left\vert \beta_{n}(u)\right\vert}{\sqrt n} <\underline{\tau}=\min
_{u\in\mathcal{L}}\tau(u)
\]
for all $n$ large enough, we only deal with$\ \rho_{+}$. Assuming that $b_{+}=1$ and $\rho_{+}(x)=xL_{+}(x)$ we have, for some $\varepsilon>0$ such that $\mathcal{L}_{\varepsilon}\subset(0,1)$ is an $\varepsilon$-neighborhood of $\mathcal{L}$, 
\begin{align*}
& \left\vert \sqrt{n}\int_{\mathcal{L}}\left(  \rho_{+}(\tau(u)+\tau
_{n}(u))-\rho_{+}(\tau(u))\right)  du-\sqrt{n}\int_{\underline{u}}%
^{\overline{u}}\rho_{+}^{\prime}(\tau(u))\tau_{n}(u)du\right\vert \\
& \leqslant\frac{\sqrt{n}}{2}\sup_{u\in\mathcal{L}_{\varepsilon}}\left\vert
\rho_{+}^{\prime\prime}(\tau(u)\right\vert \int_{\underline{u}}^{\overline{u}%
}\tau_{n}^{2}(u)du\leqslant\frac{K}{\sqrt{n}}\int_{\underline{u}}%
^{\overline{u}}\beta_{n}^{2}(u)du
\end{align*}
which almost surely vanishes by the law of the iterated logarithm. Thus we can conclude as in \cite{BFK17} by combining this with the previous Steps 1, 2, 3. In particular, the limiting variance is finite as a consequence of (\ref{CFGD_derivee}).\\
In order to complete the proof of Theorem \ref{D=R} note that whenever $F>G$ we similarly get 
\begin{align*}
& \left\vert \sqrt{n}\int_{\mathcal{L}}\left(  \rho_{-}(-\tau(u)-\tau
_{n}(u))-\rho_{-}(-\tau(u))\right)  du-\sqrt{n}\int_{\underline{u}}%
^{\overline{u}}\rho_{-}^{\prime}(-\tau(u))\tau_{n}(u)du\right\vert \\
& \leqslant\frac{K}{\sqrt{n}}\int_{\underline{u}}%
^{\overline{u}}\beta_{n}^{2}(u)du
\end{align*}
which explains why the term $\rho_{-}^{\prime}(-\tau(u))=\vert\rho_{c}^{\prime}(\tau(u))\vert$ shows up.
\subsection{The general case}

We now prove Theorem \ref{ED}. Recall that $(FG0)$ implies the existence of $0=u_{0}<u_{1}<...<u_{\kappa}=1$
such that $F^{-1}(u_{k})=G^{-1}(u_{k})$ and $A_{k}=(u_{k-1},u_{k})\subset E$
or $A_{k}\subset D$ for $k=1,...,\kappa$. We now study the mixed case where at
least one of these intervals is included in $E$ and one in $D$, so that
$\kappa\geqslant2$. Consider, using notation (\ref{IA}),
\[
\sqrt{n}(W_{c}(\mathbb{F}_{n},\mathbb{G}_{n})-W_{c}(F,G))=\sqrt{n}\sum
_{k=1}^{\kappa}I_{A_{k}}.
\]
Let $0\leqslant\lambda<\min_{1\leqslant k\leqslant\kappa}(u_{k}-u_{k-1})/2$. Define the intervals $A_{k,\lambda}^{+}=\left(  u_{k-1},u_{k-1}+\lambda\right)  \subset A_{k}$ for $2\leqslant k\leqslant\kappa$ and $A_{k,\lambda}^{-}=\left(  u_{k}-\lambda,u_{k}\right)  \subset A_{k}$ for
$1\leqslant k\leqslant\kappa-1$. If $A_{k}\subset D$ we have $F^{-1}(u)\neq
G^{-1}(u)$ for $u\in A_{k,\lambda}^{+}\cup A_{k,\lambda}^{-}$. If
$A_{k}\subset E$ the intervals $A_{k,\lambda}^{+}$ and $A_{k,\lambda}^{-}$ are
assumed to be empty instead. Consider first the intervals $A_{k,\lambda}^{+}$
for $2\leqslant k\leqslant\kappa$ and set $0<u_{-}<u_{1}<u_{\kappa-1}<u_{+}%
<1$. Since%
\[
\lim_{n\rightarrow+\infty}\sup_{u_{-}<u<u_{+}}\left|\frac{\beta_{n}(u)}{\sqrt{n}%
}\right|=0\quad a.s.
\]
we have, by $(C2)$, for $K=\sup_{u_{-}<u<u_{+}}(\rho_{-}^{\prime}(2\left\vert
\tau(u)\right\vert ),\rho_{+}^{\prime}(2\left\vert \tau(u)\right\vert
))<+\infty$,%
\[
\lim_{n\rightarrow+\infty}\sup_{u_{-}<u<u_{+}}\frac{\rho_{c}(\beta
_{n}(u)/\sqrt{n})}{|\beta_{n}(u)/\sqrt{n}|}\leqslant K\quad a.s.
\]
Therefore, in view of Step 4 in the previous proof for $F\neq G$ we get%
\[
\sqrt{n}\left\vert I_{A_{k,\lambda}^{+}}\right\vert \leqslant K\int
_{A_{k,\lambda}^{+}}\left\vert \beta_{n}(u)\right\vert du\leqslant\frac
{K}{\underline{h}}\left(  \int_{A_{k,\lambda}^{+}}\left\vert B_{n}%
(u)\right\vert du+\int_{A_{k,\lambda}^{+}}\left\vert Z_{n}(u)\right\vert
du\right)
\]
where $\underline{h}=\min_{u_{-}<u<u_{+}}\min(h_{X}(u),h_{X}(u))>0$. Lemma
\ref{Lem_E_approxim} further yields%
\[
\lim_{n\rightarrow+\infty}\mathbb{P}\left(  \sqrt{n}\left\vert I_{A_{k,\lambda
}^{+}}\right\vert >\alpha\right)  \leqslant\mathbb{P}\left(  \int
_{A_{k,\lambda}^{+}}\left\vert B(u)\right\vert du>\frac{2\alpha\underline{h}%
}{K}\right)
\]
for any $\alpha>0$ and all $2\leqslant k\leqslant\kappa$, where $B$ has the same law as $B_n=B_n^X-B_n^Y$. The latter upper
bound vanishes as $\lambda\rightarrow0$. A similar conclusion holds for
$A_{k,\lambda}^{-}$ and $1\leqslant k\leqslant\kappa-1$. Write $A_{1,\lambda
}^{\ast}=A_{1}\backslash A_{1,\lambda}^{-}$, $A_{\kappa,\lambda}^{\ast
}=A_{\kappa}\backslash A_{\kappa,\lambda}^{+}$ and $A_{k,\lambda}^{\ast}%
=A_{k}\backslash(A_{k,\lambda}^{+}\cup A_{k,\lambda}^{-})$ for $2\leqslant
k\leqslant\kappa-1$.$\medskip$

\noindent\textbf{(i)} Consider the case $1<b<2$. Fix $\lambda>0$ arbitrarily
small and write%

\begin{equation}
\sqrt{n}(W_{c}(\mathbb{F}_{n},\mathbb{G}_{n})-W_{c}(F,G))=\sqrt{n}I_{E}%
+\sqrt{n}I_{D,\lambda}^{\ast}+\sqrt{n}I_{D,\lambda}^{\pm}\label{I_etoile}%
\end{equation}
where%
\[
I_{E}=\sum_{A_{k}\subset E}I_{A_{k,\lambda}^{\ast}}=\sum_{A_{k}\subset
E}I_{A_{k}},\quad I_{D,\lambda}^{\ast}=\sum_{A_{k}\subset D}I_{A_{k,\lambda
}^{\ast}},\quad I_{D,\lambda}^{\pm}=\sum_{A_{k}\subset D}I_{A_{k,\lambda}%
^{+}\cup A_{k,\lambda}^{-}}.
\]
We just proved that%
\[
\lim_{\lambda\rightarrow0}\lim_{n\rightarrow+\infty}\mathbb{P}\left(  \sqrt
{n}I_{D,\lambda}^{\pm}>\alpha\right)  =0.
\]
Since $b>1$ we have $v_{n}/\sqrt{n}\rightarrow0$ as $n\rightarrow+\infty$.
Therefore Steps 1 to 4 of Section \ref{proofF=G} when $F=G$ show that%
\[
\lim_{n\rightarrow+\infty}\sqrt{n}I_{E}=\lim_{n\rightarrow+\infty}\frac
{\sqrt{n}}{v_{n}}v_{n}I_{E}=0\quad\text{in probability.}%
\]
In the case $\kappa\geqslant3$ then for all $2\leqslant k\leqslant\kappa-1$
with $A_{k}\subset D$ we have $\delta_{k}=\inf_{u\in A_{k}^{\ast}}\left\vert
\tau(u)\right\vert >\delta>0$ and $\tau(u)$ has constant sign on $A_{k}$. It
follows from Steps 1 to 4 of Section \ref{proofFneG} when $F\neq G$ that the weak limit
of $\sqrt{n}I_{D,\lambda}^{\ast}$ is $\int_{D_{\lambda}}\rho_{c}^{\prime}
(\tau(u))\mathbb{B}(u)du$ where $D_{\lambda}=\bigcup_{A_{k}\subset D}
I_{A_{k,\lambda}^{\ast}}$ and $\mathbb{B}(u)=B^X(u)/h_X(u)-B^Y(u)/h_Y(u)$. 
By letting $\lambda\rightarrow0$ we conclude that
\[
\sqrt{n}(W_{c}(\mathbb{F}_{n},\mathbb{G}_{n})-W_{c}(F,G))\rightarrow
_{weak}\int_{D}\rho_{c}^{\prime}(\tau(u))\mathbb{B}(u)du
\]
which is easily seen to have the normal distribution $\mathcal{N}%
(0,\sigma_{D}^{2})$.$\medskip$

\noindent\textbf{(ii)} Assume that $b=1$. Starting again from (\ref{I_etoile})
we again obtain that%
\[
\sqrt{n}I_{D,\lambda}^{\ast}\rightarrow_{weak}\int_{D_{\lambda}}\rho
_{c}^{\prime}(\tau(u))\mathbb{B}(u)du
\]
while the Steps 1 to 4 of Section \ref{proofFneG} now entails, for $v_{n}$ from
(\ref{vn}),%
\[
v_{n}I_{E}\rightarrow_{weak}\pi_{+}\int_{E}1_{\left\{  \mathbb{B}(u)>0\right\}  }%
\left\vert \mathbb{B}(u)\right\vert du+\pi_{-}\int_{E}1_{\left\{
\mathbb{B}(u)<0\right\} } \left\vert\mathbb{B}(u)\right\vert du.
\]
The above approximation with the same $\mathbb{B}$ proves that the weak convergence of the couple  $\left(\sqrt{n}I_{D,\lambda}^{\ast}, v_{n}I_{E}\right)$ holds, thus the sum weakly converges.\\

\noi Finally observe that $(C4)$ implies $\sqrt{n}/v_{n}\rightarrow L_{+}%
(0)/\pi_{+}$ and $\sqrt{n}/v_{n}\rightarrow L_{-}(0)/\pi_{-}$ as
$n\rightarrow+\infty$. As previoulsy we conclude by letting $\lambda
\rightarrow0$.

\subsection{A first special case : $F=G$ and $b=2$}
We establish Theorem \ref{W2}.
\smallskip 

\noindent\textbf{Step 0.} 
Assume $(C0)$, $\rho_{c}(x)=x^{2}$ for $\left\vert x\right\vert <x_{0}$,
$E=\mathbb{R}$, $(FG1)$, $(FG2)$ and%
\[
\lim_{u\rightarrow0}\frac{u}{h(u)}=\lim_{u\rightarrow1}\frac{1-u}%
{h(u)}=0,\quad\int_{0}^{1}\frac{u(1-u)}{h^{2}(u)}du<+\infty.
\]
This proof partially follows the line of the proof of Lemma 2.4 of \cite{Delbarrio05}.
\smallskip 

\noindent\textbf{Step 1.} 
We show that $\sup_{1/n\leqslant u\leqslant1-1/n}\left\vert \mathbb{F}%
_{n}^{-1}(u)-\mathbb{G}_{n}^{-1}(u)\right\vert \rightarrow0$ in probability,
so that the behaviour of $\rho_{c}$ near $0$ only matters. Write $h=h_{X}$.
Define $U_{i}=F(X_{i})$ and $V_{i}=F(Y_{i})$, $i=1,...,n$.
Consider $nI_{\mathcal{I}_{n}}$ with $i_{n}=1$ and%
\begin{align*}
I_{\mathcal{I}_{n}}  & =\int_{1-1/n}^{1}\left(  \mathbb{F}_{n}^{-1}%
(u)-\mathbb{G}_{n}^{-1}(u)\right)  ^{2}du\\
& =\int_{1-1/n}^{1}\left(  F^{-1}(U_{(n)})-F^{-1}(V_{(n)})\right)  ^{2}\\
& \leqslant\frac{2}{n}\left(  F^{-1}(U_{(n)})-F^{-1}(1-\frac{1}{n})\right)
^{2}du\\
& +\frac{2}{n}\left(  F^{-1}(V_{(n)})-F^{-1}(1-\frac{1}{n})\right)  ^{2}.
\end{align*}
By the mean theorem, for some random $U_{(n)}^{\ast}$ between $U_{(n)}$ and
$1-1/n$,%
\[
F^{-1}(U_{(n)})-F^{-1}(1-\frac{1}{n})=\frac{U_{(n)}-1+1/n}{h(U_{(n)}^{\ast}%
)}=\frac{U_{(n)}-1+1/n}{h(U_{(n)})}\frac{h(U_{(n)})}{h(U_{(n)}^{\ast})}.
\]
By a classical argument -- see \cite{BFK17} -- we have, thanks to $(FG2)$,
\[
\max\left(  \frac{h(U_{(n)})}{h(U_{(n)}^{\ast})},\frac{h(U_{(n)}^{\ast}%
)}{h(U_{(n)})}\right)  \leqslant\max\left(  \frac{1-U_{(n)}}{1-U_{(n)}^{\ast}%
},\frac{1-U_{(n)}^{\ast}}{1-U_{(n)}}\right)  ^{K}.
\]
Now recall that $U_{(n)}-1+1/n=O_{P}(1/n)$ and $d_{(n)}=n(1-U_{(n)}%
)\rightarrow_{weak}d_{(\infty)}$ where $d_{(\infty)}$ is a positive fintite
$r.v.$ Hence%
\begin{align*}
\left(  F^{-1}(U_{(n)})-F^{-1}(1-\frac{1}{n})\right)  ^{2}  & \leqslant
\frac{(U_{(n)}-1+1/n)^{2}}{h^{2}(U_{(n)})}\max\left(  \frac{1}{d_{(n)}%
},d_{(n)}\right)  ^{2K}\\
& =\frac{(1-U_{(n)})^{2}}{h^{2}(U_{(n)})}\left(  1-\frac{1}{d_{(n)}}\right)
^{2}\max\left(  \frac{1}{d_{(n)}},d_{(n)}\right)  ^{2K}%
\end{align*}
where $(1-U_{(n)})^{2}/h^{2}(U_{(n)})\rightarrow0$ almost surely and $\left(
1-1/d_{(n)}\right)  ^{2}\max\left(  1/d_{(n)},d_{(n)}\right)  ^{2K}=O_{P}(1)$.
Hence $nI_{\mathcal{I}_{n}}=o_{P}(1)$.
\smallskip 

\noindent\textbf{Step 2.} Now consider, for $j_{n}=n^{\beta}$,
\[
nI_{\mathcal{J}_{n}}=\int_{1-j_{n}/n}^{1-1/n}\left(  \beta_{n}^{X}%
(u)-\beta_{n}^{Y}(u)\right)  ^{2}du.
\]

\begin{lemma}
\label{Lem_b2_step2}There exists a sequence of processes $\mathbb{B}_{n}^{X}$
having the same law as 
$\mathbb{B}^{X}$ of (\ref{gauss}) such that%
\[
\Xi_{n}=\sup_{1/n\leqslant u\leqslant1-1/n}\left\vert \beta_{n}^{X}%
(u)-\mathbb{B}_{n}^{X}(u)\right\vert \frac{h(u)}{\sqrt{1-u}}=O_{P}(1)\text{.}%
\]

\end{lemma}

\textit{Proof.} It is an immediate extension of Corollary 4.2.1. page 382 of
\cite{CH93} starting from (4.2.2) of Theorem 4.2.1 of \cite{CH93}.$\quad\square$

As a consequence,
\begin{align*}
\mathbb{P}\left(  nI_{\mathcal{J}_{n}}>3\alpha\right)    & \leqslant
\mathbb{P}\left(  \int_{1-j_{n}/n}^{1-1/n}\left(  \mathbb{B}_{n}%
^{X}(u)-\mathbb{B}_{n}^{Y}(u)\right)  ^{2}du>\alpha\right)  \\
& +2\mathbb{P}\left(  \Xi_{n}^{2}\int_{1-j_{n}/n}^{1-1/n}\frac{1-u}{h^{2}%
(u)}du>\alpha\right)
\end{align*}
hence $nI_{\mathcal{J}_{n}}\rightarrow0$ in probability. We conclude the proof
by applying the Steps 3 to 5 in Section \ref{proofF=G} with many simplifications since
$L(x)=1$ now.

\subsection{A second special case : $F=G$ has compact support}

The proof of Corollary \ref{compact} follows exactly the same path as the proof of Theorem \ref{E=R} up to the following slight changes.

\noindent\textbf{Step 0.} We mainly require $(FG2)$, $(FG3)$ to apply the Hungarian construction but not $(C3)$ for the cost at $+\infty$ since the support is bounded.

\noindent\textbf{Step 1.} In Step 1 of Section \ref{proofF=G} we only need $K_n \to +\infty$.

\noindent\textbf{Step 2.} It is much shortened thanks to the boundedness of $F^{-1}$ by taking $K_n$ such that $i_n/\log\log n \to +\infty$ and (\ref{epsiln}) is no more required since by (\ref{eqcompact}) $(\sqrt{u(1-u)}/h(u))^{b'}$ is integrable.

\noindent\textbf{Steps 3 and 4.} Since $F^{-1}$ is bounded we use (\ref{eqcompact}) that implies the $a.s.$ finiteness of $\displaystyle \int_{0}^{1} \left\vert \mathbb{B}^X(u)\right\vert ^{b'} du $ and $\displaystyle \int_{0}^{1} \left\vert \mathbb{B}^Y(u)\right\vert ^{b'} du $.

\bibliographystyle{plain}
\bibliography{TCLbib_HAL}

\end{document}